\documentclass{article}
\usepackage{array}
\usepackage{amsmath}
\usepackage{amssymb}
\usepackage{amsthm}
\usepackage{bbm}
\usepackage{cancel}
\usepackage{amsfonts}
\usepackage{graphicx}
\usepackage{wrapfig}
\usepackage{caption}
\usepackage{url}
\usepackage{maplestd2e}
\theoremstyle{definition}

\DefineParaStyle{Maple Heading 1}
\DefineParaStyle{Maple Text Output}
\DefineParaStyle{Maple Dash Item}
\DefineParaStyle{Maple Bullet Item}
\DefineParaStyle{Maple Normal}
\DefineParaStyle{Maple Heading 4}
\DefineParaStyle{Maple Heading 3}
\DefineParaStyle{Maple Heading 2}
\DefineParaStyle{Maple Warning}
\DefineParaStyle{Maple Title}
\DefineParaStyle{Maple Error}
\DefineCharStyle{Maple Hyperlink}
\DefineCharStyle{Maple 2D Math}
\DefineCharStyle{Maple Maple Input}
\DefineCharStyle{Maple 2D Output}
\DefineCharStyle{Maple 2D Input}

\usepackage[top=.9in,bottom=.9in]{geometry}
 \numberwithin{equation}{section}
 
\DefineParaStyle{Maple Output}
\DefineCharStyle{2D Comment}
\DefineCharStyle{2D Math}
\DefineCharStyle{2D Output}
\begin{document}
% doesn't seem to be recognized.  \numberwithin{equation}{section}
% nor this \addtoreset{equation}{section}
%\pagestyle{empty}
% activate the following statement to obtain visible labels for each equation
%\let\labeldefs\iftrue
\let\labeldefs\iffalse
\pagestyle{plain}
%\pagenumbering{arabic}
%\newcommand{\ignore}[1]{}
\numberwithin{equation}{section}
%\widowpenalty=600
%\clubpenalty=600   !!! why doesn't "title" work???
\begin{flushleft} \vskip 0.3 in  
\centerline{\bf Further Exploration of Riemann's Functional Equation} \vskip .3in
\vskip .2in
\centerline{ Michael Milgram\footnote{mike@geometrics-unlimited.com}}

\centerline{Consulting Physicist, Geometrics Unlimited, Ltd.}
\centerline{Box 1484, Deep River, Ont. Canada. K0J 1P0}
\centerline{}
\centerline{Sept. 9, 2016}
\centerline{}
\centerline{ revised: Nov. 9, 2016 - clarifications, corrections and additions}
\centerline{}
\centerline{}
\vskip .3in
\centerline{\bf Abstract}\vskip .1in

A previous exploration of the Riemann functional equation that focussed on the critical line, is extended over the complex plane. Significant results include a simpler derivation of the fundamental equation obtained previously, and its generalization  from the critical line to the complex plane. A simpler statement of the relationship that exists between the real and imaginary components of $\zeta(s)$ and $\zeta^{\prime}(s)$ on opposing sides of the critical line is developed, reducing to a simpler statement of the same result on the critical line. An analytic expression is obtained for the sum of the arguments of $\zeta(s)$ on symmetrically opposite sides of the critical line, reducing to the analytic expression for $arg(\zeta(1/2+i\rho))$ first obtained in the previous work. Relationships are obtained between various combinations of $|\zeta(s)|$ and $|\zeta^{\prime}(s)|$, particularly on the critical line, and it is demonstrated that the difference function $arg(\zeta(1/2+i\rho))-arg(\zeta^{\prime}(1/2+i\rho))$ uniquely defines $|\zeta(1/2+i\rho)|$. A comment is made about the utility of such results as they might apply to putative proofs of Riemann's Hypothesis (RH). 

\vskip .1in

\section{Introduction}\label{sec:Intro}
In a previous report \cite[Milgram]{Milgram_Exploring}, hereinafter referred to as I, a variant of the Riemann functional equation was studied and a number of results were discovered that were either new, or well-buried in the literature. Notable was the derivation and/or discovery of:
\begin{itemize}
\item{an analytic expression for the argument of $\zeta(s)$ on the critical line $s=1/2+i\rho$ through the use of a differential equation;}
\item{a singular linear transformation that exists between the real and imaginary components of $\zeta(1/2+i\rho)$ and the corresponding components of its derivative $\zeta^{\prime}(1/2+i\rho)$;}
\item{``anomalous zeros" whose existence calls into question several well-accepted results;} and
\item{various estimates for the location and density of zeros on the critical line.}
\end{itemize}

The purpose of this work is to report on additional properties that have been found through further study of this functional equation, primarily a simplified form of both the linear transformation referred to and its derivation, its extension over the entire complex plane, and the derivation of an analytic expression for (the sum of) the argument of $\zeta(s)$ on symmetrically opposite sides of the critical line over the complex plane, again obtained through the use of a differential equation. Additionally, various relationships between $|\zeta(1/2+i\rho)|$, $|\zeta^{\prime}(1/2+i\rho)|$ as well as the real and imaginary components of $\zeta(s)$ and $\zeta^{\prime}(s)$ on opposite sides of the critical line are developed, and some criteria are deduced that must be satisfied at $s=s_0$ (where $\zeta(s_0)=0$) anywhere in the complex plane. 

\section{Recap and Notation}

The notation defined in I will be used here (see Appendix) with some extensions - particularly the specification of functional dependence. In preference to studying the complex function of a complex variable in reference to Riemann's function $\zeta(s)$ over the complex $s$ plane, I prefer to utilize its real and imaginary components, each treated as (semi-) independent real functions of a complex variable $s=\sigma+i\rho$ $(\sigma,\rho \in \Re)$, for the simple reason that many of the properties being studied must (at a fundamental computational level) be manipulated in terms of these functions. Thus, I write
\begin{equation}
\zeta(s)\equiv \zeta_R(\sigma+i\rho)+i\zeta_{I}(\sigma+i\rho)
\label{zdef}
\end{equation}
\labeldefs{\begin{flushright}$zdef$\end{flushright}}\fi
and, with due regard to the property that
\begin{equation}
\zeta(1-s)=\zeta_R(1-\sigma-i\rho)+i\zeta_{I}(1-\sigma-i\rho)=\zeta_R(1-\sigma+i\rho)-i\zeta_{I}(1-\sigma+i\rho)
\label{cmplxC}
\end{equation}
\labeldefs{\begin{flushright}$cmplxC$\end{flushright}}\fi
and to distinguish results valid over the entire complex $s$ plane from those valid only on the critical line $\sigma=1/2$, I write
$\zeta(s)$ to mean the former and $\zeta$ to mean the latter, with the extension that, for any (relevant) function,
\begin{equation}
\widetilde \zeta(s)\equiv\zeta(1-\sigma+i\rho)\neq\zeta(1-\sigma-i\rho),
\label{Zdef}
\end{equation}
\labeldefs{\begin{flushright} $ZdefC$\end{flushright}}\fi
and $\rho>0$ always. Throughout, derivatives denoted by $^{\prime}$ refer to the operation $\frac{\partial}{\partial\rho}$ unless the argument of the operand is specified as $s$, in which case it refers to $\frac{d}{ds}$.\newline

Thus by considering its real and imaginary parts, (see the Appendix for a summary of symbols), the traditional form of Riemann's functional equation reads

\begin{eqnarray}
{\it \widetilde {\zeta}_R}(s)=\frac{1}{2\,\pi_\sigma}(\,{ {g_{{2}}(s)\zeta_{{R}}(s)}}+{{g_{{1}}(s)\zeta_{{I}}(s)}}) \nonumber \\
{\it \widetilde {\zeta}_I}(s)=\frac{1}{2\,\pi_\sigma}(\,{ {g_{{1}}(s)\zeta_{{R}}}(s)}-{{g_{{2}}(s)\zeta_{{I}}(s)}})
\label{basic}
\end{eqnarray}
\labeldefs{\begin{flushright}$basic$\end{flushright}}\fi

or, more succinctly
\begin{eqnarray}
{\it \widetilde {\zeta}_R}(s)=\frac{1}{2\,\pi_\sigma}\zeta_p(s) \\
{\it \widetilde {\zeta}_I}(s)=-\frac{1}{2\,\pi_\sigma}\zeta_m(s)
\label{Simple}
\end{eqnarray}
\labeldefs{\begin{flushright}$Simple$\end{flushright}}\fi

from the polar form of which (see Appendix and I) the relationship between the arguments of $\zeta(s)$ on symmetrically opposite sides of the critical line immediately follows: 
\begin{equation}
\mapleinline{inert}{2d}{tan(alpha_a) = -(tan(alpha)*g[2]-g[1])/(tan(alpha)*g[1]+g[2])}
{\[\displaystyle \tan \left( {\it \widetilde{\alpha}(s)} \right)\equiv\,{\it \widetilde {\zeta}_I}(s)/{\it \widetilde {\zeta}_R}(s)  =-{\frac {\tan \left( \alpha(s) \right) g_{{2}}(s)-g_{{1}}(s)}{\tan \left( \alpha(s) \right) g_{{1}}(s)+g_{{2}}(s)
\mbox{}}}\]} \,.
\label{tans}
\end{equation}
\labeldefs{\begin{flushright}tans\end{flushright}}\fi

% for the ffollowing, see bottom of file invert_2.mw

For completeness' sake, the inverse of \eqref{basic} is
\begin{eqnarray}
\mapleinline{inert}{2d}{Zeta[R] =
1/8*Pi[sigma]*(Zeta[a,I]*g[1](s)+Zeta[a,R]*g[2](s))/abs(GAMMA(s))^2/c[
0];}{%
\[
{\zeta _{R}(s)}={\displaystyle \frac {1}{8}} \,{\displaystyle 
\frac {{\pi _{\sigma }}\,({\widetilde{\zeta}_{\,I}(s)}\,{g_{1}}(s) + {\widetilde{\zeta} 
_{\,R}(s)}\,{g_{2}}(s))}{ \left|  \! \,\Gamma (s)\, \!  \right| 
^{2}\,{c_{0}}}} \label{basic_R}
\]
}
\\
\mapleinline{inert}{2d}{Zeta[I] =
-1/8*Pi[sigma]*(Zeta[a,I]*g[2](s)-Zeta[a,R]*g[1](s))/abs(GAMMA(s))^2/c
[0];}{%
\[
{\zeta _{I}(s)}= - {\displaystyle \frac {1}{8}} \,{\displaystyle 
\frac {{\pi _{\sigma }}\,({\widetilde{\zeta}_{\,I}(s)}\,{g_{2}}(s) - {\widetilde{\zeta} 
_{\,R}}(s)\,{g_{1}}(s))}{ \left|  \! \,\Gamma (s)\, \!  \right| 
^{2}\,{c_{0}}}\,.} 
\]
}
\label{basic_I}
\end{eqnarray}
\labeldefs{\begin{flushright}$basic_R,basic_I$\end{flushright}}\fi

% from file absZratio.mw%
Squaring \eqref{basic_R} and \eqref{basic_I} then adding, eventually produces an expression equivalent to relatively well-known (see for example \cite[Spira, Eq.(2)]{Spira1965}) expressions for the ratio of magnitudes of $\zeta(s)$ on opposite sides of the critical line:
\begin{equation}
\mapleinline{inert}{2d}{abs(Zeta)^2/abs(ZetaA)^2 = (1/2)*(2*Pi)^(2*sigma)/((cos(Pi*sigma)+cosh(Pi*rho))*abs(GAMMA(s))^2)}{\[\displaystyle {\frac { \displaystyle | \zeta(s) |^{2}}{  | \widetilde{\zeta}(s)|^{2}}}={\frac { \left( 2\,\pi \right) ^{2\,\sigma}}{2\, \left( \cos \left( \pi\,\sigma \right) +\cosh \left( \pi\,\rho \right)  \right) 
\mbox{}  \left| \Gamma \left( s \right)  \right|^{2}}}\]}\equiv\Phi(s)\,.
\label{Zmags}
\end{equation}
\labeldefs{\begin{flushright}Zmags\end{flushright}}\fi
For $\sigma=1/2$ in \eqref{Zmags}, $\Phi=1$ (see \cite[NIST, Eq. (5.4.4)] {NIST}).\newline

The variant form of the functional equation, as utilized in I is

\begin{equation}
\mapleinline{inert}{2d}{Zeta(1,1-s)/Zeta(1,s)+2*cos(1/2*Pi*s)*GAMMA(s)*(2*Pi)^(-s) =
f(s)*Zeta(1-s)/Zeta(1,s);}{%
\[
\mathfrak{L}(s)\equiv{\displaystyle \frac {\zeta^{\prime} (1 - s)}{\zeta^{\prime} (s)}}  + \chi(s)
={\displaystyle \frac {\mathit{f}(s)\,\zeta (1 - 
s)}{\zeta^{\prime} (s)}}\equiv\mathfrak{T}(s) 
\]
}
\label{variant}
\end{equation}
\labeldefs{\begin{flushright}variant\end{flushright}}\fi

where
\begin{equation}
\mapleinline{inert}{2d}{chi(s) = 2*cos(1/2*Pi*s)*GAMMA(s)/((2*Pi)^s);}{%
\[
\chi (s)={\displaystyle \frac {2\,\mathrm{cos}({\displaystyle 
\frac {\pi \,s}{2}} )\,\Gamma (s)}{(2\,\pi )^{s}}} 
\]
}
\label{chi}
\end{equation}
\labeldefs{\begin{flushright}chi\end{flushright}}\fi

valid for all $s$. Recall that in I it was proven that $\zeta^{\prime}(1/2+i\rho)\neq0$, except possibly at a zero, and, assuming the Riemann Hypothesis(RH), Spira \cite{Spira1973} has shown that $\zeta^{\prime}(s)\neq 0$ for all $\sigma<1/2$. 
\newline

Throughout, in an attempt to reduce results with many terms into a comprehensible whole, I adhere to the convention that any symbol containing one of the 8 primitive $\zeta$ functions (real and imaginary components of $\zeta(s),\widetilde{\zeta}(s)$ and derivatives thereof) somewhere in its structure will always be represented by a variation of the letter $\zeta$, whereas if a symbol does not contain the letter $\zeta$, it is usually a function of other variables, notably $\Gamma$ and its derivatives, as well as trigonometric and hyperbolic functions of the variables $\sigma$ and $\rho$, with the notable exception that the argument(s) of $\zeta$ and/or $\zeta^{\prime}$ may also appear. As well, many of the calculations are rather lengthy and require the use of a computer algebra program. Here I use the computer program Maple \cite{Maple} extensively and include the annotation ``(Maple)" at the appropriate location(s) as the justification and source of a particular result.  

\section{Over the whole plane}

\subsection{The real and imaginary parts of $\zeta(s)$}

Many of the following results require fairly lengthy derivation and considerable manipulation using a computer algebra program. By equating the real and imaginary parts of \eqref{variant}, it is possible to relate the real and imaginary functions $\widetilde{\zeta}_R(s)$ and $\widetilde{\zeta}_I(s)$ with the other 6 components - see \eqref{ZAR0} and \eqref{ZAI0}. After incorporating \eqref{basic}, we find (Maple) a number of interesting, useful and simpler variations: 

\begin{equation}
\mapleinline{inert}{2d}{Zeta[aI] =
4*((-g[1]*p[1]+g[2]*p[2])*Zeta[1,I]+(-g[1]*p[2]-g[2]*p[1])*Zeta[1,R])*
c[0]/((2*Pi)^sigma)/(p[1]^2+p[2]^2)-8*(Zeta[1,aI]*p[2]+Zeta[1,aR]*p[1]
)*c[0]/(p[1]^2+p[2]^2);}{%
\[
\mathit{\widetilde{\zeta}_R(s)}=-{\displaystyle \frac{4\,c_0}{({p_{1}}^{2}(s) + {p_{2}}^{2}(s))}}\left( {\displaystyle \frac{1}{(2\,\pi )^{
\sigma }}}{\displaystyle  {[h_2(s)\,{\zeta^{\prime}_{I}(s)} + h_1(s)\,{\zeta^{\prime}_{R}(s)}]\,{}}{}}  - {\displaystyle {
2\,[{\widetilde{\zeta}^{\prime} _{{I}}(s)}\,{p_{1}(s)} - {\widetilde{\zeta}^{\prime}_{{
R}}(s)}\,{p_{2}(s)}]\,{}}{}} \right)
\]
}
\label{ZAR}
\end{equation}

\labeldefs{\begin{flushright}ZAR\end{flushright}}\fi
and

\begin{equation}
\mapleinline{inert}{2d}{Zeta[aI] =
4*((-g[1]*p[1]+g[2]*p[2])*Zeta[1,I]+(-g[1]*p[2]-g[2]*p[1])*Zeta[1,R])*
c[0]/((2*Pi)^sigma)/(p[1]^2+p[2]^2)-8*(Zeta[1,aI]*p[2]+Zeta[1,aR]*p[1]
)*c[0]/(p[1]^2+p[2]^2);}{%
\[
\mathit{\widetilde{\zeta}_I(s)}={\displaystyle \frac{4\,c_0}{({p_{1}(s)}^{2} + {p_{2}(s)}^{2})}}\left( {\displaystyle \frac{1}{(2\,\pi )^{
\sigma }}}{\displaystyle  {[h_1\,{\zeta^{\prime}_{I}(s)} - h_2(s)\,{\zeta^{\prime}_{R}(s)}]\,{}}{}}  - {\displaystyle {
2\,[{\widetilde{\zeta}^{\prime} _{{I}}(s)}\,{p_{2}(s)} + {\widetilde{\zeta}^{\prime}_{{
R}}(s)}\,{p_{1}(s)}]\,{}}{}} \right)
\]
}
\label{ZAI}
\end{equation}
\labeldefs{\begin{flushright}ZAI\end{flushright}}\fi
%see EqGen5.mw Eq(26) %
which together give the relationship between the components of $\widetilde{\zeta}(s)$ and their counterparts symmetrically across the critical line, demonstrating, as suggested in \eqref{variant}, that knowledge of the derivatives on both sides of the critical line, specifies $\widetilde{\zeta}(s)$ and hence $\zeta(s)$ itself via the transformation $\sigma\rightarrow1-\sigma$. In I, it was noted that when relationships such as \eqref{ZAR} and \eqref{ZAI} are limited to the critical line, the transformation relationship is singular and thus non-invertible (e.g. see \eqref{ZR} and \eqref{ZI} below). For arbitrary values of $s\neq 1/2+i\rho$ however, we find (Maple) the inverted transformation that defines $\zeta^{\prime}(s)$ in terms of $\widetilde{\zeta}(s)$ and $\widetilde{\zeta}^{\prime}(s)$ - its components on the opposite side of the critical line:

%see EqGen5.mw Eq(30) %

\begin{equation}
\mapleinline{inert}{2d}{T2a1 := Zeta[1,R] =
-1/64*(8*CC*Zeta[1,a,I]*g[1]+8*CC*Zeta[1,a,R]*g[2]+Zeta[a,I]*g[1]*p[2]
+Zeta[a,I]*g[2]*p[1]+h[1](s)*Zeta[a,R])/abs(GAMMA)^2/CC^2*TPis;}{%
\maplemultiline{
{\zeta^{\prime}_{\,R}(s)}= 
 - {\displaystyle \frac {(2\pi)^\sigma}{64\, \left|  \! \,\Gamma \, \!  \right| ^{2
}\,\mathit{c_0^2}}} \,{\displaystyle  {\left( 8\,
\mathit{c_0}\,[\,{\widetilde{\zeta}^{\prime}_{\,I}(s)}\,{g_{1}(s)} + {
\widetilde{\zeta}^{\prime}_{\,R}(s)}\,{g_{2}(s)}] + {h_{2}(s)}\,{\widetilde{\zeta}_{\,I}(s)}\, + {h_{1}}(s)\,{\widetilde{\zeta} _{
\,R}(s)}\right) }{}}  }
}
\label{Z1R}
\end{equation}
\labeldefs{\begin{flushright}Z1R\end{flushright}}\fi
and
\begin{equation}
\mapleinline{inert}{2d}{T1a1 := Zeta[1,I] =
1/64*(8*CC*Zeta[1,a,I]*g[2]-8*CC*Zeta[1,a,R]*g[1]-h[2](s)*Zeta[a,R]+h[
1](s)*Zeta[a,I])/abs(GAMMA)^2/CC^2*TPis;}{%
\[
{\zeta^{\prime}_{I}(s)}={\displaystyle \frac {(2\,\pi)^\sigma}{64\, \left|  \! \,\Gamma \, \!  \right| ^{2}\,\mathit{c_0^2}}} 
\,\left( {\displaystyle  {8\,\mathit{c_0}\,[\,{\widetilde{\zeta}^{\prime}_{I}(s)}\,
{g_{2}(s)} - {\widetilde{\zeta}^{\prime}_{R}(s)}\,{g_{1}(s)}] - {h_{2
}}(s)\,{\widetilde{\zeta}_{R}(s)} + {h_{1}}(s)\,{\widetilde{\zeta}_{I}(s)}\,, 
}}\,\right )\, 
\]
}
\label{Z1I}
\end{equation}
\labeldefs{\begin{flushright}Z1I\end{flushright}}\fi
along with the inverse(s):
%see file Gen_squared3.mw%
\begin{equation}
\mapleinline{inert}{2d}{Z1aR = -(1/2)*(Zeta__I(1, s)*g[1](s)+Zeta__R(1, s)*g[2](s))/Pi__sigma-(1/16)*(-h[3](s)*Zeta__I(s)+h[4](s)*Zeta__R(s))/(Pi__sigma*c__0)}{\[\displaystyle {\widetilde{\zeta}_R^{\prime}(s)}=-{\frac { g_{{1}} \left( s \right)\zeta^{\prime}_{I} \left(s \right) +\zeta^{\prime}_{R} \left(s \right) g_{{2}} \left( s \right) }{2\,\pi_{\sigma
\mbox{}}}}+{\frac {h_{{3}} \left( s \right) \zeta_{I} \left( s \right) -h_{{4}} \left( s \right) \zeta_{R} \left( s \right) }{16\,\pi_{\sigma
\mbox{}}\,c_{0}}}\]}
\label{Za1R}
\end{equation}
\labeldefs{\begin{flushright}Za1R\end{flushright}}\fi
and
\begin{equation}
\mapleinline{inert}{2d}{Z1aI = (1/2)*(Zeta__I(1, s)*g[2](s)-g[1](s)*Zeta__R(1, s))/Pi__sigma+(1/16)*(h[4](s)*Zeta__I(s)+h[3](s)*Zeta__R(s))/(Pi__sigma*c__0)}{\[\displaystyle {\widetilde{\zeta}_I^{\prime}(s)}={\frac { g_{{2}} \left( s \right)\zeta^{\prime}_{I} \left(s \right) -g_{{1}} \left( s \right) \zeta^{\prime}_{R} \left( s \right) }{2\,\pi_{\sigma
\mbox{}}}}+{\frac {h_{{4}} \left( s \right) \zeta_{I} \left( s \right) +h_{{3}} \left( s \right) \zeta_{R} \left( s \right) }{16\,\pi_{\sigma
\mbox{}}\,c_{0}}}\]}\,.
\label{Za1I}
\end{equation}
\labeldefs{\begin{flushright}Za1I\end{flushright}}\fi

Further, by appropriate (and patient) manipulation of the same equations, we find (Maple)

\begin{align}\nonumber
\left|  \! \,\zeta^{\prime} (s)\, \!  \right| ^{2}=\; &
{\displaystyle \frac{\pi _{\sigma }}{g_{
1}(s)}} 
 \displaystyle  {\left( [{\zeta^{\prime} _{\,I}(s)}\,{\widetilde{\zeta} _{\,R}(s)} + {
\zeta^{\prime}_{\,R}(s)}\,{\widetilde{\zeta}_{\,I}(s)}]\,{q_{1}(s)} + [ - {\zeta^{\prime}_{
\,I}(s)}\,{\widetilde{\zeta}_{\,I}(s)} + {\zeta^{\prime}_{\,R}(s)}\,{\widetilde{\zeta} _{\,R}(s)}]
\,{q_{2}(s)} \right. }  \\ 
 & \left. - 2\,{\zeta^{\prime}_{\,I}(s)}\,{\widetilde{\zeta}^{\prime}_{\,R}(s)}  - 2\,{
\zeta^{\prime}_{\,R}(s)}\,{\widetilde{\zeta}^{\prime}_{\,I}(s)}\,\right)\,.  
\label{Z1abs1}
\end{align}
\labeldefs{\begin{flushright}Z1abs1\end{flushright}}\fi

% see file Eqgen_5.mw%
Applying the polar form of the the various elements to the ratio of \eqref{ZAR} and \eqref{ZAI}, we find (Maple) an equivalent form of \eqref{tans}
\begin{equation}
\mapleinline{inert}{2d}{tan(alpha[a,s]) =
(exp(-1/2*Pi*rho)*sin(-1/2*Pi*sigma-alpha(s)+rho[theta(s)])+exp(1/2*Pi
*rho)*sin(1/2*Pi*sigma-alpha(s)+rho[theta(s)]))/(exp(1/2*Pi*rho)*cos(1
/2*Pi*sigma-alpha(s)+rho[theta(s)])+exp(-1/2*Pi*rho)*cos(-1/2*Pi*sigma
-alpha(s)+rho[theta(s)]));}{%
\[
\mathrm{tan}({\widetilde{\alpha}(s)})={\displaystyle \frac {e^{ - 
\frac {\pi \,\rho }{2}}\,\mathrm{sin}( - {\displaystyle  {
\pi \,\sigma }/{2}}  - \alpha (s) + {\rho _{\theta (s)}}) + e^{
\frac {\pi \,\rho }{2}}\,\mathrm{sin}({\displaystyle  {\pi 
\,\sigma }/{2}}  - \alpha (s) + {\rho _{\theta (s)}})}{e^{\frac {
\pi \,\rho }{2}}\,\mathrm{cos}({\displaystyle  {\pi \,
\sigma }/{2}}  - \alpha (s) + {\rho _{\theta (s)}}) + e^{ - 
\frac {\pi \,\rho }{2}}\,\mathrm{cos}( - {\displaystyle  {
\pi \,\sigma }/{2}}  - \alpha (s) + {\rho _{\theta (s)}})}\,.} 
\]
}
\label{tan_alpha}
\end{equation}
\labeldefs{\begin{flushright}$tan_{alpha}$\end{flushright}}\fi
In a sense, \eqref{tans} and/or \eqref{tan_alpha} are functional equations for the argument of $\zeta(s)$.
\newline

% see lambda in file "zeros_off_the_critical_line.mw"
Similarly, applying the polar forms to \eqref{Z1R} and \eqref{Z1I} we find
\begin{equation}
\mapleinline{inert}{2d}{tan(beta) =
-(8*(tan(beta[1])*g[2](s)-g[1](s))*lambda*c[0]+h[1](s)*tan(alpha[1])
-h[2](s))/((8*tan(beta[1])*g[1](s)+8*g[2](s))*lambda*c[0]+h[2](s)*tan(
alpha[1])+h[1](s));}{%
\[
\mathrm{tan}(\beta(s))= - {\displaystyle \frac {8\,\widetilde{\zeta}_{r}(s) \,{c_{0}}[\,\mathrm{tan}({
\widetilde{\beta}(s)})\,{g_{2}}(s) - {g_{1}}(s)] + {h
_{1}}(s)\,\mathrm{tan}({\widetilde{\alpha}(s)}) - {h_{2}}(s)}{8\,\widetilde{\zeta}_{r}(s) \,{c_{0}[\,\mathrm{
tan}({\widetilde{\beta}(s)})\,{g_{1}}(s) + {g_{2}}(s)]\,
} + {h_{2}}(s)\,\mathrm{tan}({\widetilde{\alpha}(s)}) + {h_{1}}(s)}} 
\]
}
\label{tanB}
\end{equation}
\labeldefs{\begin{flushright}tanB\end{flushright}}\fi
and the inverse
\begin{equation}
\mapleinline{inert}{2d}{tan(beta(s)[1]) = -((8*tan(beta(s))*g[2](s)-8*g[1](s))*c__0*Zeta__r(s)+tan(alpha(s))*h[4](s)+h[3](s))/((8*tan(beta(s))*g[1](s)+8*g[2](s))*c__0*Zeta__r(s)-tan(alpha(s))*h[3](s)+h[4](s))}{\[\displaystyle \tan(\widetilde{\beta} \left( s \right)) =-{\frac {  8\,c_{0}\,\zeta_{r} \left( s \right)\left( \tan \left( \beta \left( s \right)  \right) g_{{2}} \left( s \right) -g_{{1}} \left( s \right)  \right) 
\mbox{}+\tan \left( \alpha \left( s \right)  \right) h_{{4}} \left( s \right) +h_{{3}} \left( s \right) }{8\,c_{0}\,\zeta_{r} \left( s \right) \left( \tan \left( \beta \left( s \right)  \right) g_{{1}} \left( s \right) +g_{{2}} \left( s \right)  \right)  
\mbox{}-\tan \left( \alpha \left( s \right)  \right) h_{{3}} \left( s \right) +h_{{4}} \left( s \right) }}\]}
\label{tanB1}
\end{equation}
\labeldefs{\begin{flushright}tanB1\end{flushright}}\fi

where
\begin{equation}
\widetilde{\zeta}_{r}(s)\equiv\widetilde{\zeta}^{\prime}_R(s)/\widetilde{\zeta}_R(s)
\label{zetaRR_def}
\end{equation}
\labeldefs{\begin{flushright}$zetaRR_def$\end{flushright}}\fi 
and
\begin{equation}
{\zeta}_{r}(s)\equiv{\zeta}^{\prime}_R(s)/{\zeta}_R(s).
\label{zetaR_def}
\end{equation}
\labeldefs{\begin{flushright}$zetaR_def$\end{flushright}}\fi 
% see file Gen_Squared3.mw%

Solving \eqref{tanB} gives
\begin{equation}
\mapleinline{inert}{2d}{lambda =
((-tan(alpha[1](s))-tan(beta(s)))*h[1](s)+(-tan(alpha[1](s))*tan(beta(
s))+1)*h[2](s))/((8*tan(beta[1](s))*tan(beta(s))-8)*g[1](s)+(8*tan(bet
a[1](s))+8*tan(beta(s)))*g[2](s))/c[0];}{%
\[
\widetilde{\zeta}_{r}(s) ={\displaystyle \frac {( - \mathrm{tan}({\widetilde{\alpha}}(s))
 - \mathrm{tan}(\beta (s)))\,{h_{1}}(s) + ( - \mathrm{tan}({
\widetilde{\alpha}}(s))\,\mathrm{tan}(\beta (s)) + 1)\,{h_{2}}(s)}{8\,c_{0}\left((
\mathrm{tan}({\widetilde{\beta}}(s))\,\mathrm{tan}(\beta(s)) - 1)\,{g_{
1}}(s) + (\mathrm{tan}({\widetilde{\beta}}(s)) + \mathrm{tan}(
\beta (s)))\,{g_{2}}(s)\right)\,{}}} 
\]
}
\label{ZetaRR}
\end{equation}
\labeldefs{\begin{flushright}ZetaRR\end{flushright}}\fi
and/or solving \eqref{tanB1} gives
\begin{equation}
\mapleinline{inert}{2d}{Z__r = (1/8)*((tan(alpha(s))*h[3](s)-h[4](s))*tan(beta[1](s))-tan(alpha(s))*h[4](s)-h[3](s))/(c__0*((tan(beta[1](s))*g[1](s)+g[2](s))*tan(beta(s))+g[2](s)*tan(\widetilde{beta}(s))-g[1](s)))}{\[\displaystyle \zeta_{r}(s)={\frac { \left( \tan \left( \alpha \left( s \right)  \right) h_{{3}} \left( s \right) -h_{{4}} \left( s \right)  \right) \tan ( \widetilde{\beta}\left( s  \right) -\tan \left( \alpha \left( s \right)  \right) h_{{4}} \left( s \right) -h_{{3}} \left( s \right) }{8\,c_{0}\, \left(  \left( \tan ( \widetilde{\beta}\left( s )  \right) g_{{1}} \left( s \right) +g_{{2}} \left( s \right)  \right) \tan \left( \beta \left( s \right)  \right) +g_{{2}} \left( s \right) \tan ( \widetilde{\beta} \left( s )  \right) -g_{{1}} \left( s \right)  \right) 
\mbox{}}}\]}
\label{Zeta_rr}
\end{equation}
\labeldefs{\begin{flushright}$Zeta_rr$\end{flushright}}\fi

results that may be useful elsewhere.

\subsection{The argument of $\zeta(s)$} 
\begin{figure}[h] % from file "alpha_p_graphics.mw"%
\centering
\includegraphics[width=.6\textwidth]{{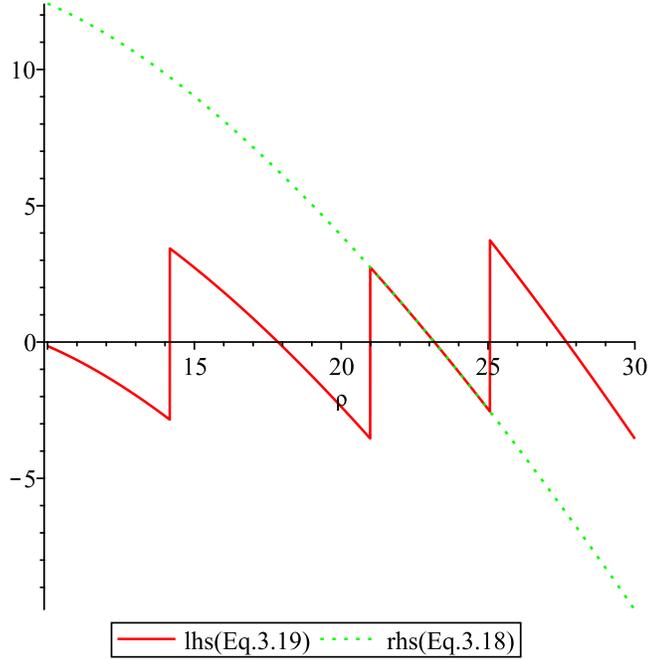}}
\caption{Plot of the left and right-hand sides of \eqref{argmean} and \eqref{argmean2} using $\sigma=1/3$ and $k=2$.}
\label{fig:alpha_p}
\end{figure}

%see file Gen6.mw%
Motivated by the derivation of Eq. (6.1) of I (see \eqref{Zmain} below), multiply together \eqref{basic_R},\eqref{Z1R},\eqref{basic_I} and \eqref{Z1I} in corresponding pairs to find (Maple)
\begin{equation}
\mapleinline{inert}{2d}{Z1I*ZI+Z1R*ZR = -(Z1aI*ZaI+Z1aR*ZaR)*abs(Zeta)^2/abs(Zeta_a)^2-(1/8)*p[2]*abs(Zeta)^2/c[0]}{\[\displaystyle {\zeta^{\prime}_I(s)}\,{\zeta_I(s)}+{\zeta^{\prime}_R(s)}\,{\zeta_R(s)}=-{\frac { \left( \displaystyle {\widetilde{\zeta^{\prime}}_I(s)}\,{\widetilde{\zeta}_I(s)}+{\widetilde{\zeta^{\prime}}_R(s)}\,{\widetilde{\zeta}_R(s)} \right)   \left| \zeta(s) \right|  ^{2}
\mbox{}}{ \left| { \widetilde{\zeta}}(s) \right|  ^{2}}}
\mbox{}-\,{\frac {p_{{2}}  \left| \zeta(s) \right|  ^{2}}{8\,c_{{0}}}}\]}\,,
\label{arg1}
\end{equation}
\labeldefs{\begin{flushright}arg1\end{flushright}}\fi 
which, after substitution (see Appendix) and simplifications (including \eqref{Zmags}), can be written more symmetrically as
\begin{equation}
\mapleinline{inert}{2d}{Zc[2]/abs(Zeta_a)^2+Zc[1]/abs(Zeta)^2 = -(1/4)*(2*Psi[2]*c[0]-Pi*sin(Pi*sigma))/c[0]}{\[\displaystyle {\frac {\displaystyle {\widetilde{\zeta^{\prime}}_I(s)}\,{\widetilde{\zeta}_I(s)}+{\widetilde{\zeta^{\prime}}_R(s)}\,{\widetilde{\zeta}_R(s)}}{ | \widetilde{\zeta}(s) |  ^{2}}}+{\frac {\displaystyle {\zeta^{\prime}_I(s)}\,{\zeta_I(s)}+{\zeta^{\prime}_R(s)}\,{\zeta_R(s)}}{  \left| \zeta(s) \right| ^{2}}}=-{\frac {2\,\Psi_{{2}}c_{{0}}
\mbox{}-\pi\,\sin \left( \pi\,\sigma \right) }{4\,c_{{0}}}}\]}\,,
\label{deq2}
\end{equation}
\labeldefs{\begin{flushright}deq2\end{flushright}}\fi 
being the generalization of Eq. (6.1) of I over the complex plane. Employing the same logic carefully presented in Section 6 of I, and taking the derivatives with respect to $\rho$, it is evident that \eqref{deq2} can be interpreted as a differential equation for the sum of the (continuous) arguments $\alpha_p(s)\equiv\alpha(s)+\,\widetilde{\alpha}(s)$, specifically
\begin{equation}
\mapleinline{inert}{2d}{Diff(alpha[p](\sigma+i\,rho), rho) = ln(2*Pi)-Re(Psi(i*rho+sigma))+(1/2)*Pi*sin(Pi*sigma)/(cos(Pi*sigma)+cosh(Pi*rho))}{\[\displaystyle {\frac {\partial}{{\partial}\rho}}\alpha_{{p}} \left(\sigma+i\, \rho \right) =\ln  \left( 2\,\pi \right) -\Re \left( \psi \left( \sigma+i\rho\right)  \right) +\frac{\pi}{2}\,{\frac {\sin \left( \pi\,\sigma \right) }{\cos \left( \pi\,\sigma \right) +\cosh \left( \pi\,\rho \right) }}\]}
\label{deq}
\end{equation}
\labeldefs{\begin{flushright}deq\end{flushright}}\fi 

whose solution (Maple) gives

\begin{align}
\displaystyle {\alpha_{p}(s)}= 
-\Re\int_{0}^{\rho}\!\psi \left( \sigma+it \right) \,{\rm d}t 
  -\arctan \left( {\frac { \left( \cos \left( \pi\,\sigma \right) -1 \right) \sinh \left( \pi\,\rho/2 \right) }{\cosh \left( \pi\,\rho/2 \right) 
\mbox{}\sin \left( \pi\,\sigma \right) }} \right) +\rho_{\pi} +k\pi
\label{argmean}
\end{align}
\labeldefs{\begin{flushright}argmean\end{flushright}}\fi 
and, for a specific value of $k$,
\begin{equation}
arg(\zeta(s))+arg(\widetilde{\zeta}(s))=\alpha_{p}(s),
\label{argmean2}
\end{equation}
\labeldefs{\begin{flushright}argmean2\end{flushright}}\fi 

thereby expressing the sum of $arg(\zeta(s))$ on opposite sides of the critical line, in terms of simple, basic, well-known functions, reducing to Eq. (6.9) of I when $\sigma=1/2$, taking into account the identity {\it $$ \arctan(e^{\rho\pi})-\arctan (\tanh (\rho \pi/2))=\pi/4\,.$$} See Figure \ref{fig:alpha_p} for an example, and note that the right-hand side of \eqref{deq} is always negative for $\rho\gg0$, implying that, for reasonably large values of $\rho$, the locus of $\zeta(s)$ always travels in a clockwise direction with increasing $\rho$ in the complex $\zeta(s)$ plane (see Section 9.1 of I).\newline

In the presence of a discontinuity, the value of $k$ included in \eqref{argmean} is effectively a winding number in the complex $\zeta(s)$ plane - each time either $arg(\zeta(s))$ or $arg(\widetilde{\zeta}(s))$ changes from $-\pi$ to $+\pi$ (with increasing $\rho$) as the locus of $\zeta(s)$ crosses the negative axis ($\zeta_R(s)<0$ or $\widetilde{\zeta}_R(s)<0$), $k$ will increase by 2 in order to maintain (the veracity of) \eqref{argmean2} over a limited (but continuous) range of $\rho$. If RH is true, or complex zeros (if such exist when $\sigma\neq1/2$) are always of even order, odd values of $k$ will never occur, because at a discontinuity associated with a simple (or odd-order) zero, $arg(\zeta(s))$ only jumps by odd multiples of $\pi$ (see Eq. (8.2) of I). Thus, if all zeros of $\zeta(s)$ are simple, the presence of an odd value of $k$ in \eqref{argmean} {\bf over a continuous range of $\rho$ (as opposed to a numerical solution at a single value of $\rho$)} when $\sigma\neq1/2$ would correspond to a counter-example to the truth of RH, in which case \eqref{argmean} could possibly be utilized as the basis for a test for RH. \newline

Furthermore, although the two results \eqref{tan_alpha} and \eqref{argmean} do not supply a solution for $\alpha(s)$ and $\widetilde{\alpha}(s)$ individually, after applying the identity \eqref{tanId} (see below), they do yield the following identity:
%see file gen6.mw%
\begin{align}\label{ident}
& \displaystyle \Im \left(\rm {Log\Gamma} \left( s \right)  \right) -\arctan \left(\tan\left(\frac{\pi\sigma}{2} \right )\,\tanh\left( \frac{\pi\rho}{2}\right) \right) -\rho\,\ln  \left( 2\,\pi \right) -k\pi \\
& =\arctan \left( {\frac {\sin \left( 1/2\,\pi\,\sigma \right) \cos \left( \rho_{{\theta \left( s \right) }} \right) \sinh \left( 1/2\,\pi\,\rho \right) 
\mbox{}+\cos \left( 1/2\,\pi\,\sigma \right) \sin \left( \rho_{{\theta \left( s \right) }} \right) \cosh \left( 1/2\,\pi\,\rho \right) }{\sin \left( 1/2\,\pi\,\sigma \right) \sin \left( \rho_{{\theta \left( s \right) }} \right) \sinh \left( 1/2\,\pi\,\rho \right) 
\mbox{}-\cos \left( 1/2\,\pi\,\sigma \right) \cos \left( \rho_{{\theta \left( s \right) }} \right) \cosh \left( 1/2\,\pi\,\rho \right) }} \right)  \\
& =-\arctan(\frac{g_1(s)}{g_2(s)})
\end{align}
\labeldefs{\begin{flushright}ident\end{flushright}}\fi 

which in turn reduces to the trivial identity \eqref{argGammaS} if the $arctan$ addition rule \cite[NIST, Eq.(4.24.15)]{NIST} is applied, where (see Eq. (6.12) of I)
\begin{equation}
\mapleinline{inert}{2d}{Re(int(Psi(sigma+I*t), t = 0 .. rho)) = Im(lnGAMMA(sigma+I*rho))}{\[\displaystyle \Re \left( \int_{0}^{\rho}\!\psi \left( \sigma+it \right) \,{\rm d}t \right) =\Im \left( {\rm Log\Gamma}
\mbox{} \left( \sigma+i\rho \right)  \right) \]} \,.
\label{IntPsi}
\end{equation}
\labeldefs{\begin{flushright}IntPsi\end{flushright}}\fi 

%file lnGamma_extensions.mw%
In a similar vein, note that the left-hand side of \eqref{argmean} is invariant under the substitution $\sigma\rightarrow 1-\sigma$ and the right-hand side is not. Thus, after performing that substitution and subtracting, we find
\begin{equation}
\mapleinline{inert}{2d}{tan(Im(lnGAMMA(sigma+I*rho))-Im(lnGAMMA(1-sigma+I*rho))) = -tanh(Pi*rho)/tan(Pi*sigma)}{\[\displaystyle \tan \left( \Im \left[ {\rm Log\Gamma} \left( \sigma+i\rho \right)  -{\rm Log\Gamma} \left( 1-\sigma+i\rho \right)  \right] 
\mbox{} \right) =-{\frac {\tanh \left( \pi\,\rho \right) }{\tan \left( \pi\,\sigma \right) }}\]},
\label{tanIm}
\end{equation}
\labeldefs{\begin{flushright}tanIm\end{flushright}}\fi

or, equivalently

\begin{equation}
\cos(\theta(s)-\widetilde{\theta}(s))=\frac{\sin(\pi\sigma)}{|\cos(\pi\sigma)|\sqrt{\tan^2(\pi\sigma)+\tanh^2(\pi\rho)}}
\label{cosTd}
\end{equation}
\labeldefs{\begin{flushright}cosTd\end{flushright}}\fi 

neither of which appear in the usual references \cite[NIST, Section 5]{NIST}. In the limit $\rho\gg1$, \eqref{cosTd} becomes

\begin{equation}
\cos(\theta(s)-\widetilde{\theta}(s))\sim{\sin(\pi\sigma)}\,.
\end{equation}

Operating on \eqref{tanIm} with $\frac{\partial}{\partial\sigma}$ gives  

\begin{equation}
\mapleinline{inert}{2d}{Im(Psi(sigma+I*rho)+Psi(1-sigma+I*rho)) = -Pi*sinh(2*Pi*rho)/(-cosh(2*Pi*rho)+cos(2*Pi*sigma))}{\[\displaystyle \Im \left[ \psi \left( \sigma+i\rho \right) +\psi \left( 1-\sigma+i\rho \right)  \right] ={\frac {\pi\,\sinh \left( 2\,\pi\,\rho \right) 
\mbox{}}{\cosh \left( 2\,\pi\,\rho \right) -\cos \left( 2\,\pi\,\sigma \right) }}\]}
\label{ImPsi}
\end{equation}
\labeldefs{\begin{flushright}ImPsi\end{flushright}}\fi

reducing to a known result \cite[NIST, Eq. (5.4.17)] {NIST} when $\sigma=1/2$. Similarly, operating on \eqref{tanIm} with $\frac{\partial}{\partial\rho}$ gives

\begin{equation}
\mapleinline{inert}{2d}{Im(Psi(sigma+I*rho)+Psi(1-sigma+I*rho)) = -Pi*sinh(2*Pi*rho)/(-cosh(2*Pi*rho)+cos(2*Pi*sigma))}{\[\displaystyle \Re \left[ \psi \left( \sigma+i\rho \right) -\psi \left( 1-\sigma+i\rho \right)  \right] ={\frac {-\pi\,\sin \left( 2\,\pi\,\sigma \right) 
\mbox{}}{\cosh \left( 2\,\pi\,\rho \right) -\cos \left( 2\,\pi\,\sigma \right) }}\]}\,.
\label{RePsi}
\end{equation}
\labeldefs{\begin{flushright}RePsi\end{flushright}}\fi
Although neither of these appear in the usual references, combining \eqref{ImPsi} and \eqref{RePsi} yields the standard reflection formula \cite[NIST, Eq. (5.5.4)]{NIST} for $\psi(\sigma+i\rho)$, and therefore both can alternatively be derived by working backwards from that relationship. See also Srinivasan and Zvengrowski \cite{Srin&Zven}.

\section{On the Critical Line}
\subsection{Relationships devolving from $\alpha(s)$}

In the case that $\sigma=1/2$, we have $\widetilde{\alpha}(s)=\alpha(s)$, and solving \eqref{tan_alpha} gives (Maple)

\begin{equation}
\mapleinline{inert}{2d}{tan(alpha) = (exp(Pi*rho)*sin((1/4)*Pi+rho[theta])-cos((1/4)*Pi+rho[theta]))/(exp(Pi*rho)*cos((1/4)*Pi+rho[theta])+sqrt(exp(2*Pi*rho)+1)+sin((1/4)*Pi+rho[theta]))}{\[\displaystyle \tan \left( \alpha \right) ={\frac {{{\rm e}^{\pi\,\rho}}\sin \left( \pi/4+\rho_{{\theta}} \right) -\cos \left( \pi/4+\rho_{{\theta}} \right) 
\mbox{}}{{{\rm e}^{\pi\,\rho}}\cos \left( \pi/4+\rho_{{\theta}} \right) + \sqrt{{{\rm e}^{2\,\pi\,\rho}}+1}+\sin \left( \pi/4+\rho_{{\theta}} \right) }}\]}
\label{tana2}
\end{equation}
\labeldefs{\begin{flushright}$tana2$\end{flushright}}\fi
or equivalently
\begin{equation}
\mapleinline{inert}{2d}{tan(alpha) = -(exp(Pi*rho)*cos((1/4)*Pi+rho[theta])-sqrt(exp(2*Pi*rho)+1)+sin((1/4)*Pi+rho[theta]))/(exp(Pi*rho)*sin((1/4)*Pi+rho[theta])-cos((1/4)*Pi+rho[theta]))}{\[\displaystyle \tan \left( \alpha \right) =-{\frac {{{\rm e}^{\pi\,\rho}}\cos \left( \pi/4+\rho_{{\theta}} \right) - \sqrt{{{\rm e}^{2\,\pi\,\rho}}+1}+\sin \left( \pi/4+\rho_{{\theta}} \right) 
\mbox{}}{{{\rm e}^{\pi\,\rho}}\sin \left( \pi/4+\rho_{{\theta}} \right) -\cos \left( \pi/4+\rho_{{\theta}} \right) }}\]}\,.
\label{tana1}
\end{equation}
\labeldefs{\begin{flushright}$tana1$\end{flushright}}\fi

On the critical line, the analytic representation of $\alpha$ and therefore $arg(\zeta)$ (see \eqref{argmean} or Eq. (6.15) of I) is
\begin{equation}
\mapleinline{inert}{2d}{argument(Zeta(1/2+t*I)) = -1/2*int(Re(psi(1/2+rho*I)),rho = 0 ..
t)+1/2*t*ln(2)+1/2*t*ln(Pi)-9/8*Pi+1/2*arctan(exp(Pi*t));}{%
\maplemultiline{
\alpha=
  - {\displaystyle \frac {1}{2}} \,{\displaystyle \int _{0}^{\rho}} 
\Re (\psi ({\displaystyle  {1}/{2}}  + i\,t))\,dt  + 
{\displaystyle \frac {\rho}{2}} \,\mathrm{ln}(2\,\pi)  - 
{\displaystyle \frac {9\,\pi }{8}}  + {\displaystyle \frac {1}{2}
} \,\mathrm{arctan}(e^{\pi \,\rho})+k\,\pi \,.}
}
\label{ArgZt=0}
\end{equation}
\labeldefs{\begin{flushright}ArgZt=0\end{flushright}}\fi{.}
\newline

so, for example, \eqref{tana1} can also be written 
\begin{equation}
\mapleinline{inert}{2d}{tan(1/2*rho[theta]-1/8*Pi+1/2*arctan(exp(Pi*rho))) =
(exp(Pi*rho)*cos(1/4*Pi+rho[theta])+sin(1/4*Pi+rho[theta])-(exp(2*Pi*r
ho)+1)^(1/2))/(-sin(1/4*Pi+rho[theta])*exp(Pi*rho)+cos(1/4*Pi+rho[thet
a]));}{%
\[
\mathrm{tan}({\displaystyle \frac {1}{2}} \,{\rho _{\theta }} - 
{\displaystyle \frac {\pi }{8}}  + {\displaystyle \frac {1}{2}} 
\,\mathrm{arctan}(e^{\pi \,\rho}))={\displaystyle \frac {e^{
\pi \,\rho }\,\mathrm{cos}({\displaystyle  {\pi }/{4}}  + {
\rho _{\theta }}) + \mathrm{sin}({\displaystyle  {\pi }/{4}} 
 + {\rho _{\theta }}) - \sqrt{e^{2\,\pi \,\rho} + 1}}{ - 
\,e^{\pi \,\rho }\,\mathrm{sin}({\displaystyle {\pi }/{4}}  + {\rho _{\theta }}
) + \mathrm{cos}({\displaystyle {\pi }/{4
}}  + {\rho _{\theta }})}\,.} 
\]
}
\label{tanEq}
\end{equation}
\labeldefs{\begin{flushright}tanEq\end{flushright}}\fi
Because of the definition \eqref{rhoTheta}, both \eqref{ident} and \eqref{tanEq} can be interpreted as functional equations for any of $\theta(\rho),\theta(s)$ or the imaginary part of the $\mathrm {LogGamma}$ function respectively (see \eqref{IntPsi}), without reference to the $\zeta$ function at all. In contrast to \eqref{ArgZt=0} with \eqref{IntPsi}, and based on a well-known result (reproduced as Eq. (2.10) of I), Brent has recently \cite[Section 4]{2016arXiv160903682B} obtained
\begin{equation}
\mapleinline{inert}{2d}{-argument(Zeta(1/2+I*rho)) = (1/2)*theta+(1/2)*rho*ln(2*Pi)-(1/8)*Pi+arctan(exp(-Pi*rho))}{\[\displaystyle {-arg} \left( \zeta \left( 1/2+i\rho \right)  \right) =-arg(\Gamma(1/2+i\rho))/2+\frac{\rho}{2}\,\ln  \left( 2\,\pi \right) 
\mbox{}+\pi/8-\frac{1}{2}\,\arctan \left( {{\rm e}^{-\pi\,\rho}} \right)\]}\,,
\label{Brent}
\end{equation}
\labeldefs{\begin{flushright}Brent\end{flushright}}\fi
omitting a term equal to $-k\pi$. As discussed in I, the term $k\pi$ in \eqref{ArgZt=0} tallies the zeros of $\zeta(1/2+i\rho)$. This would also be true of \eqref{Brent} only if the term $arg(\Gamma(1/2+i\rho))$ in \eqref{Brent} were to be interpreted as the continuous function ${\rm \Im(Log\Gamma(1/2+i\rho)}$) (see \cite[comment following Eq. (1)]{2016arXiv160903682B}). Figure \ref{fig:Brent_Fig} illustrates the difference between the two interpretations (also see Appendix, comment following \eqref{argZeta2}).
\newline

\begin{figure}[h] % from file "Brent&Carma.mw"%
\centering
\includegraphics[width=.6\textwidth]{{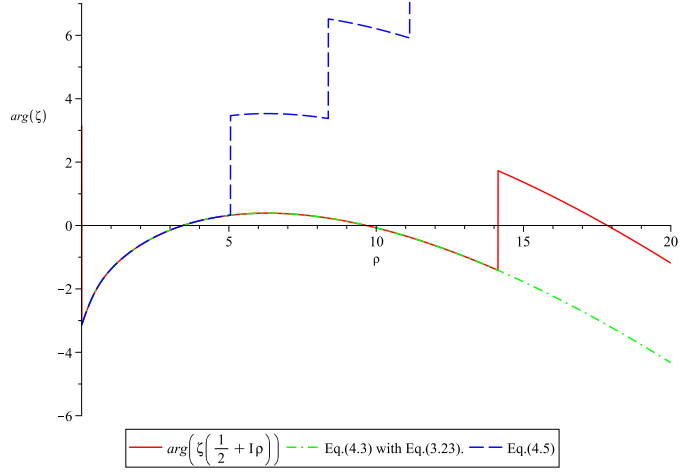}}
\caption{$arg(\zeta(1/2+i\rho))$ in the vicinity of the first zero, showing the effect of two interpretations of the first term on the right-hand side of \eqref{Brent}, as well as the necessity of including terms $\pm k\pi$. Both results in the Figure set $k=0.$}
\label{fig:Brent_Fig}
\end{figure}

Also, on the critical line $\sigma=1/2$, $\mathit{\widetilde{\zeta}_I}=\mathit{{\zeta}_I}$, $\mathit{\widetilde{\zeta}_R}=\mathit{{\zeta}_R}$ and similarly $\mathit{\widetilde{\zeta}^{\prime}_I}=\mathit{{\zeta}^{\prime}_I}$ and $\mathit{\widetilde{\zeta}^{\prime}_R}=\mathit{{\zeta}^{\prime}_R}$. Substituting these identifications into \eqref{ZAR} and \eqref{ZAI} we find
\begin{equation}
\mapleinline{inert}{2d}{Zeta[R]/f = (b+1/2)*Zeta[1,R]+a*Zeta[1,I];}{%
\[
{\displaystyle \frac {{\zeta _{R}}}{f}} =(b + {\displaystyle 
\frac {1}{2}} )\,{\zeta^{\prime} _{R}} + a\,{\zeta^{\prime} _{I}}
\]
}
\label{ZR}
\end{equation}
\labeldefs{\begin{flushright}ZR\end{flushright}}\fi 
and
\begin{equation}
\mapleinline{inert}{2d}{Zeta[I]/f = a*Zeta[1,R]-(b-1/2)*Zeta[1,I];}{%
\[
{\displaystyle \frac {{\zeta _{I}}}{f}} =a\,{\zeta^{\prime} _{R}} - (
b - {\displaystyle \frac {1}{2}} )\,{\zeta^{\prime} _{I}}
\]
}
\label{ZI}
\end{equation}
\labeldefs{\begin{flushright}ZI\end{flushright}}\fi
where $a$ and $b$ are defined in the Appendix. The results \eqref{ZR} and \eqref{ZI} are of the same, but simpler (and equivalent) form compared to Eqs. (4.9) of I. As in I, these two results define a linear, singular transformation between $\zeta^{\prime}$ and $\zeta$ on the critical line, because, as is easily shown, the determinant of the transformation matrix
\begin{equation}
a^2+b^2-1/4=0.
\label{aSq-bSq}
\end{equation}
\labeldefs{\begin{flushright}aSq-bSq\end{flushright}}\fi

%see file tan(beta).mw
Further, it is possible to identify the functions $a$ and $b$ by first transforming \eqref{ZR} and \eqref{ZI} into their polar counterparts, then writing
\begin{equation}
\cos^{2}(\alpha)=\zeta^2_R/(\zeta^2_I+\zeta^2_R)
\label{cosDef}
\end{equation}
\labeldefs{\begin{flushright}cosDef\end{flushright}}\fi
and substituting the right-hand sides of  \eqref{ZR} and \eqref{ZI} along with \eqref{aSq-bSq} into \eqref{cosDef} to identify
\begin{equation}
%b=\cos^{2}(\alpha)-1/2
b=\cos(2\alpha)/2
\label{Bid}
\end{equation}
\labeldefs{\begin{flushright}Bid\end{flushright}}\fi
from which we correspondingly find
\begin{equation}
a=\sin(2\alpha)/2\,.
\label{Aid}
\end{equation}
\labeldefs{\begin{flushright}Aid\end{flushright}}\fi

With reference to \eqref{a-def} and \eqref{b-def}, and, relative to \eqref{tana1}, there exists a simpler relationship between the polar angles of $\Gamma(1/2+i\rho)$ and $\zeta(1/2+i\rho)$, that being
\begin{equation}
\mapleinline{inert}{2d}{sin(2*alpha) =
(cos(rho[theta])*sinh(1/2*Pi*rho)-sin(rho[theta])*cosh(1/2*Pi*rho))/co
sh(Pi*rho)^(1/2);}{%
\[
\mathrm{sin}(2\,\alpha )={\displaystyle \frac {\mathrm{cos}({\rho
 _{\theta }})\,\mathrm{sinh}({  {\frac{1}{2}\pi \,\rho }{
}} ) - \mathrm{sin}({\rho _{\theta }})\,\mathrm{cosh}(
{ {\frac{1}{2}\pi \,\rho }{}} )}{\sqrt{\mathrm{cosh}(\pi
 \,\rho )}}\,.} 
\]
}
\label{sinA=T}
\end{equation}
\labeldefs{\begin{flushright}sinA=T\end{flushright}}\fi

Utilizing Eq. (2.9) of I, then yields the following:% see file Lam2.mw"%
\begin{equation}
\mapleinline{inert}{2d}{cos(2*alpha) = (cosh((1/2)*Pi*rho)*cos(rho[theta])+sinh((1/2)*Pi*rho)*sin(rho[theta]))/sqrt(cosh(Pi*rho))}{\[\displaystyle \cos \left( 2\,\alpha \right) ={\frac {\cosh \left( \frac{1}{2}\,\pi \,\rho \right) \cos \left( \rho_{{\theta}} \right) +\sinh \left( \frac{1}{2}\,\pi \,\rho \right) \sin \left( \rho_{{\theta}} \right) 
\mbox{}}{ \sqrt{\cosh \left( \pi \,\rho \right) }}}\]},
\label{cosA=T}
\end{equation}
\labeldefs{\begin{flushright}cosA=T\end{flushright}}\fi
thereby reducing Eq.(2.9) of I to a trivial trigonometric identity:
\begin{equation}
\tan(\alpha)=1/\sin(2\alpha) - 1/\tan(2\alpha)\,.
\label{tanId}
\end{equation}
\labeldefs{\begin{flushright}tanId\end{flushright}}\fi
 
The set \eqref{sinA=T} and \eqref{cosA=T} can be inverted, giving
%see file lambda.mw%
\begin{equation}
\mapleinline{inert}{2d}{cos(rho[theta]) = sqrt(cosh(Pi*rho))*(sinh((1/2)*Pi*rho)*sin(2*alpha)+cosh((1/2)*Pi*rho)*cos(2*alpha))/(sinh((1/2)*Pi*rho)^2+cosh((1/2)*Pi*rho)^2)}{\[\displaystyle \cos \left( \rho_{{\theta}} \right) ={\frac {  \sinh \left( \pi\,\rho/2 \right) \sin \left( 2\,\alpha \right) +\cosh \left( \pi\,\rho/2 \right) \cos \left( 2\,\alpha \right)  
\mbox{}}{\sqrt{\cosh(\pi\rho)} }}\]}
\label{cosr1}
\end{equation}
\labeldefs{\begin{flushright}cosr1\end{flushright}}\fi
and \begin{equation}
\mapleinline{inert}{2d}{sin(rho[theta]) = sqrt(cosh(Pi*rho))*(sinh((1/2)*Pi*rho)*cos(2*alpha)-cosh((1/2)*Pi*rho)*sin(2*alpha))/(sinh((1/2)*Pi*rho)^2+cosh((1/2)*Pi*rho)^2)}{\[\displaystyle \sin \left( \rho_{{\theta}} \right) ={\frac {   \sinh \left( \pi\,\rho/2 \right) \cos \left( 2\,\alpha \right) -\cosh \left( \pi\,\rho/2 \right) \sin \left( 2\,\alpha \right)  
\mbox{}}{ \sqrt{\cosh \left( \pi\,\rho \right) }  }}\]}
\label{sinr1}
\end{equation}
\labeldefs{\begin{flushright}sinr1\end{flushright}}\fi
or, together
% see file tan(beta).mw
\begin{equation}
\mapleinline{inert}{2d}{tan(rho[theta]) =
-((exp(Pi*rho)+1)*tan(2*alpha)-exp(Pi*rho)+1)/((exp(Pi*rho)-1)*tan(2*a
lpha)+exp(Pi*rho)+1);}{%
\[
\mathrm{tan}({\rho _{\theta }})= - {\displaystyle \frac {(e^{\pi
 \,\rho } + 1)\,\mathrm{tan}(2\,\alpha ) - e^{\pi \,\rho } + 1
}{(e^{\pi \,\rho } - 1)\,\mathrm{tan}(2\,\alpha ) + e^{\pi \,
\rho } + 1}} \,
\]
}
\label{tanT=A}
\end{equation}
\labeldefs{\begin{flushright}tanT=A\end{flushright}}\fi
and, asymptotically (as $\rho\rightarrow\infty$),
\begin{equation}
\tan(\rho_{\theta})\sim \frac{1-\tan(2\alpha)}{1+\tan(2\alpha)}\,.
\label{tanAsy}
\end{equation}
\labeldefs{\begin{flushright}tanAsy\end{flushright}}\fi

% see file Eq4p17_check.mw%
As expected, \eqref{cosr1}-\eqref{tanT=A} reduce to tautologies with the identification \eqref{ArgZt=0}. Because the denominator of \eqref{tanT=A} vanishes when 
\begin{equation}
\tan(2\alpha)=-\coth(\rho\pi/2)
\label{Tan2A}
\end{equation}
\labeldefs{\begin{flushright}tan2A\end{flushright}}\fi
and the numerator doesn't, discontinuities in $\rho_{\theta}$ are related to $arg(\zeta)$ through numerical solutions of \eqref{Tan2A}. Thus, to the extent that $\coth(\rho\pi/2)\approx 1$, discontinuities in $\rho_{\theta}$ (and hence $\arg(\Gamma)$) will occur at those values of $\rho$ where $arg(\zeta)$ passes through $-\pi/8$. 
\newline

\subsection{Relationships involving $\zeta$ and its derivatives}
The fundamental relationships (I, Eqs. (6.1) and (6.6)) are reproduced here:
\begin{equation}
\zeta_R^2+\zeta_I^2=f(\zeta^{\prime}_I\zeta_I+\zeta^{\prime}_R\zeta_R)
\label{Zmain}
\end{equation}
\labeldefs{\begin{flushright}Zmain\end{flushright}}\fi 
and
\begin{equation}
\alpha^{\prime}=1/f\,.
\label{alphaP}
\end{equation}
\labeldefs{\begin{flushright}alphaP\end{flushright}}\fi 

The majority of the results of I were based on \eqref{Zmain}, which was derived in a very complicated manner \cite[see Supplemental\_Material]{Milgram_Exploring}. As an alternative to a reduction of \eqref{arg1}, \eqref{Zmain} can now be simply obtained from the new (and independent) results \eqref{ZI}, \eqref{ZR}, \eqref{Aid} and \eqref{Bid} by multiplying \eqref{ZI} by $\zeta^{\prime}_I$, \eqref{ZR} by $\zeta^{\prime}_R$, then adding the two to yield the intermediate identity
\begin{equation}
\mapleinline{inert}{2d}{Zeta[I]*Zeta1[I]+Zeta[R]*Zeta1[R] = -f*(Zeta1[I]^2-Zeta^{\prime}[R]^2)*cos(alpha)^2+2*Zeta1[R]*f*Zeta1[I]*sin(alpha)*cos(alpha)+Zeta1[I]^2*f}{\[\displaystyle \zeta_{{I}} \zeta^{\prime}_I+\zeta_{{R}}{ \zeta^{\prime}}_R=f\left[ \,- \left( {\zeta^{\prime}_I}^{2}-{\zeta^{\prime}_R}^{2} \right)   \cos^2 \left( \alpha \right)  
\mbox{}+\zeta^{\prime}_{{R}}\zeta^{\prime}_I\sin \left(2 \alpha \right) +{{\zeta^{\prime}}_I}^{2}\right]\]}\,,
\label{Zmain_int}
\end{equation}
\labeldefs{\begin{flushright}Zmain\_int\end{flushright}}\fi

converting all trigonometric factors back into the components $\zeta_R$, $\zeta_I$ and $|\zeta|$ using \eqref{ZetaExp} and factoring the resulting expression. \eqref{Zmain} is an immediate consequence.\newline

Further useful results can be easily obtained by squaring  \eqref{ZR} and \eqref{ZI}, applying \eqref{Bid} and \eqref{Aid} and simplifying, all of which eventually lead to

\begin{eqnarray}
%\zeta_R^2/f^2=\frac{|\zeta^{\prime}|^2}{4}(1+\cos(2\alpha))(1+\cos(2(\alpha-\beta)))\\%
\zeta_R^2/f^2={|\zeta^{\prime}|^2}\cos^2(\alpha)\,\cos^2(\alpha-\beta) \label{ZrSq} \\
%\zeta_I^2/f^2=\frac{|\zeta^{\prime}|^2}{4}(1-\cos(2\alpha))(1+\cos(2(\alpha-\beta)))%
\zeta_I^2/f^2={|\zeta^{\prime}|^2}\sin^2(\alpha)\cos^2(\alpha-\beta)\label{ZiSq}
\end{eqnarray}
\labeldefs{\begin{flushright}ZrSq,ZiSq\end{flushright}}\fi

and

\begin{equation}
\frac{|\zeta|^2}{|\zeta^{\prime}|^2}=f^2\,\cos^2(\alpha-\beta)\,.
\label{ZSqZpSq}
\end{equation}
\newline
\labeldefs{\begin{flushright}ZSqZpSq0,ZSqZpSq\end{flushright}}\fi

A more general form of \eqref{ZSqZpSq} can be obtained by direct differentiation of \eqref{Zeta-polar} using \eqref{alphaP}, yielding 
\begin{equation}
\frac{|\zeta|}{|\zeta^{\prime}|}=f\,\cos(\alpha-\beta)\,,
\label{ZZp}
\end{equation}
\labeldefs{\begin{flushright}ZZp\end{flushright}}\fi

the polar form of \eqref{Zmain}. However, this procedure cannot be utilized as a derivation of \eqref{Zmain} because \eqref{alphaP} is not independently known without \eqref{Zmain}. At this point, we also note that \eqref{ZetaRR} reduces (Maple) to \eqref{ZZp} on the critical line. This is a  fairly lengthy reduction that makes use of \eqref{cosr1} and \eqref{sinr1} as well as many trigonometric identities; an interim noteworthy result is
% file lambda.mw%
\begin{align}
&{\displaystyle{\frac { \left( -\sinh \left( \pi\,\rho \right) \sin \left( 2\,\rho_{{\theta}} \right) -\cos \left( 2\,\rho_{{\theta}} \right)  \right) \cos \left( 4\,\beta \right) + \left( -\cos \left( 2\,\rho_{{\theta}} \right) \sinh \left( \pi\,\rho \right) +\sin \left( 2\,\rho_{{\theta}} \right)  \right) \sin \left( 4\,\beta \right)
\mbox{}+\cosh \left( \pi\,\rho \right) }{\cosh \left( \pi\,\rho \right) \cos \left( \alpha-\beta \right) -\sinh \left( \pi\,\rho \right) \sin \left( 2\,\rho_{{\theta}}+3\,\beta+\alpha \right) -\cos \left( 2\,\rho_{{\theta}}+3\,\beta+\alpha \right) }}} \\
&=2\,\cos \left( \alpha-\beta \right)\,,\nonumber 
\label{interim}
\end{align}
\labeldefs{\begin{flushright}interim\end{flushright}}\fi
which, it should be noted, is independent of $\beta$, although not obviously so.
\newline

From \eqref{ZZp}, for $\rho\gtrsim\rho_s=6.28\dots$, where $f<0$, (see Eq.(2.5) of I), we have  
\begin{equation}
\cos(\alpha-\beta)\leq0\,.
\label{abIneq}
\end{equation}
\labeldefs{\begin{flushright}abIneq\end{flushright}}\fi

From Eq. (7.8) of I it is known that $\zeta(1/2+i\rho)=0$ iff 
\begin{equation}
\alpha-\beta=(n+1/2\,)\pi\,.
\label{al-be}
\end{equation}
\labeldefs{\begin{flushright}al-be\end{flushright}}\fi

In addition to locating the full zero via the criterion \eqref{al-be}, \eqref{ZrSq} and \eqref{ZiSq} together demonstrate that the real half-zero $\zeta_R=0,\,\zeta_I\neq0$ occurs when
\begin{equation}
\alpha=(n+1/2\,)\,\pi,
\label{alRe}
\end{equation}
\labeldefs{\begin{flushright}alRe\end{flushright}}\fi
and the imaginary half-zero $\zeta_I=0,\,\zeta_R\neq0$ occurs when
\begin{equation}
\alpha=n\pi\,, \hspace{1cm} {\rm where }\;\;\; n=0,\pm1,\dots\,\,,
\label{alIm}
\end{equation}
\labeldefs{\begin{flushright}alIm\end{flushright}}\fi
thereby generalizing Section 13 of I.  Also, recall that ``anomalous zeros" (see I) are characterized by the imaginary half-zero $\zeta_I=0, \zeta_R<0, \zeta^{\prime}_R>0$. From \eqref{ZrSq} and \eqref{ZiSq}, it is clear that, assuming that the zeros of $\zeta$ are simple, the combination $\alpha(\rho_0)=\pm\,\pi/2, \beta(\rho_0)=0$ cannot occur, since, in that case $\zeta_R\rightarrow0$ faster than $\zeta_I$.
\newline

Further interesting relationships can be obtained (\cite[see Supplemental\_Material]{Milgram_Exploring}) from the supplement to I where, following Eq.(9), we find the identity
\begin{equation}
\mapleinline{inert}{2d}{abs(Zeta)^2/(abs(Zeta1)^2*f^2)+(-Zeta[I]*Zeta1[R]+Zeta[R]*Zeta1[I])^2/(abs(Zeta1)^2*abs(Zeta)^2) = 1}{\[\displaystyle {\frac {   \left| \zeta  \right|  ^{2}}{  \left| {\it \zeta^{\prime}} \right|  ^{2}{f}^{2}}}+{\frac { \left( -{\it \zeta^{\prime}}_{{R}}\zeta _{{I}}+{\it \zeta^{\prime}}_{{I}}\zeta _{{R}} \right) ^{2}}{  \left| {\it \zeta^{\prime}} \right|   ^{2}  \left| \zeta  \right| ^{2}}}
\mbox{}=1\]}\,.
\label{Eq9a}
\end{equation}
\labeldefs{\begin{flushright}Eq9a\end{flushright}}\fi

By straightforward differentiation, the numerator of the second term can be identified as
\begin{equation}
{ \left( -{\it \zeta^{\prime}}_{{R}}\zeta _{{I}}+{\it \zeta^{\prime}}_{{I}}\zeta _{{R}} \right) ^{2}}=\left( \dfrac{d}{d\rho} |\zeta|^2 \right)^2/4=|\zeta|^2\left(|\zeta|^{\prime}\right)^2
\label{Eq9adiff}
\end{equation}
\labeldefs{\begin{flushright}Eq9adiff\end{flushright}}\fi
so that
\begin{equation}
{(|\zeta|^{\prime})^2}{}=|\zeta^{\prime}|^2-\frac{|\zeta|^2}{f^2}
\label{Zprime}
\end{equation}
\labeldefs{\begin{flushright}Zprime\end{flushright}}\fi

which can alternatively be written 
%see file "attempt_at_de.mw Eq(4)%

\begin{eqnarray}
 \left({\displaystyle \dfrac{d}{d\rho}\,\log|\zeta|}\,\right)^2=&{\displaystyle \frac{|\zeta^{\prime}|^2}{|\zeta|^2}-\frac{1}{f^2}}
\label{LogZSq}\\
=& {\displaystyle \frac{\tan^2(\alpha-\beta)}{f^2}}
\label{LogZSqA}
\end{eqnarray}
\labeldefs{\begin{flushright}LogZSq,LogZSqA\end{flushright}}\fi

after applying \eqref{ZSqZpSq}. This generalizes Eq.(11.3) of I. In another form, \eqref{LogZSqA} can be rewritten
\begin{equation}
\frac{|\zeta|^{\prime}}{|\zeta^{\prime}|}=\sin(\alpha-\beta)
\label{ZpZp}
\end{equation}
\labeldefs{\begin{flushright}ZpZp\end{flushright}}\fi
obtained in the same way as \eqref{ZZp}; by transforming to polar form (see \eqref{Zeta-polar}), \eqref{ZpZp} also reduces to a simple trigonometric identity, which also means that \eqref{LogZSqA} can be rewritten as

\begin{equation}
 {\displaystyle \dfrac{d}{d\rho}\,\log|\zeta|}= {\displaystyle \frac{\tan(\alpha-\beta)}{f}}\,.
\label{LogZA}
\end{equation}
\labeldefs{\begin{flushright}LogZA\end{flushright}}\fi

Integrating \eqref{LogZA} shows that $|\zeta|$ is defined by the difference of the arguments $\alpha$ and $\beta$, that is

\begin{equation}
\frac{|\zeta(1/2+i\rho_{2})|}{|\zeta(1/2+i\rho_{1})|}=\exp\left(\int_{\rho_{1}}^{\rho_2}\frac{\tan(\alpha-\beta)}{f}\,d\rho\right)\,.
\label{ExpRep'n}
\end{equation}
\labeldefs{\begin{flushright}ExpRep'n\end{flushright}}\fi
The result \eqref{ExpRep'n} has been numerically verified in several cases where $(\rho_1,\rho_2)$ does not encompass $\rho_0$ - also see Eq.(11.5) of I. From all the above, and using \eqref{alphaP}, we find

\begin{equation}
\mapleinline{inert}{2d}{Eqb := diff(beta(rho), rho) = 1/f-(diff(f(rho), rho))/(f*tan(alpha-beta))-(Diff(abs(Zeta1(rho)), rho))/(tan(alpha-beta)*abs(Zeta1(rho)))}{\[\displaystyle \beta^{\prime} \left( \rho \right) =\frac{2}{f}-
{\frac {f^{\prime} }{f\tan \left( \alpha-\beta \right) }}-{\frac { \left| {\it \zeta^{\prime}}   \right|^{\prime} }{\tan \left( \alpha-\beta \right)  \left| {\it \zeta^{\prime}}   \right| }}\]}
\label{diffBeta}
\end{equation}
\labeldefs{\begin{flushright}diffBeta\end{flushright}}\fi
which can be rewritten in the more intriguing form

\begin{equation}
\mapleinline{inert}{2d}{(Diff(beta(rho), rho))/(Diff(alpha(rho), rho)) = 2-(Diff(ln(abs(Zeta[1])*f), rho))/(Diff(ln(abs(Zeta)), rho))}{\[\displaystyle {\frac {\beta^{\prime} \left( \rho \right) }{\alpha^{\prime} \left( \rho \right) }}=2-{\frac {\left( \ln  \left(  \left| \zeta^{\prime} \right| f \right) \right)^{\prime} }{ \left( \ln \left| \zeta \right|  \right )^{\prime} }}\]}\,.
\label{DbDa}
\end{equation}
\labeldefs{\begin{flushright}DbDa\end{flushright}}\fi

In Eq.(4.7) of I, the following term arose

\begin{equation}
\mapleinline{inert}{2d}{L[1](rho) =
(Pi/cosh(Pi*rho))^(1/2)*((-1/2*Pi^(1/2)/cosh(Pi*rho)-psi[2]/Pi^(1/2))*
abs(Zeta[1])^2+(2*sin(beta)*zeta[2,I]+2*cos(beta)*zeta[2,R])*abs(zeta[
1])/Pi^(1/2))*T1;}{%
\maplemultiline{
{L}_{1}(\rho)= 
\,\left[ -{\displaystyle\frac{1}{ f(\rho)}} + {\displaystyle  {(\mathrm{sin}(\beta 
)\,{\zeta^{\prime\prime} _{\;I}} + \mathrm{cos}(\beta )\,{\zeta^{\prime\prime} _{\;R}}
)\, / \left|  \! \,{\zeta^{\prime}}\, \!  \right| }} \right ]\,
 }
}
\label{L1}
\end{equation}
\labeldefs{\begin{flushright}L1\end{flushright}}\fi

and it was claimed that $L_1(\rho)<0$, based on a numerical study. Here I add the claim-reinforcing observation that a sign change in $L_1$ would imply the existence of a zero of $\zeta(1/2+i\rho)$, which in turn would lead to an inconsistency between the solutions of Eqs. (4.3) and (4.8) of I. Further rearrangement using Eq.(6.2) of I and converting from polar form shows that \eqref{L1} can be rewritten as

\begin{equation}
L_1=\displaystyle { {\beta^{\prime} \left( \rho \right) }-\alpha^{\prime}(\rho)}<0\,.
\label{L1Neg}
\end{equation}
\labeldefs{\begin{flushright}L1Neg\end{flushright}}\fi

\begin{figure}[h] % from file "lnGamma_extensions.mw"%
\centering
\includegraphics[width=.6\textwidth]{{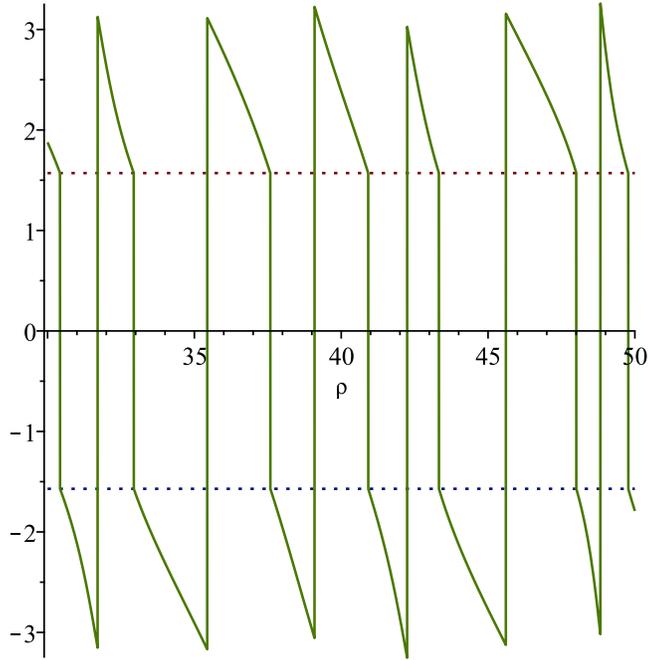}}
\caption{Plot of $\beta(\rho)-\alpha(\rho)$ over the range $30\leq\rho\leq50$. Note that the slope is always negative (Eq. \eqref{L1Neg}), there exist discontinuities corresponding to $\zeta(1/2+i\rho)=0$ at $\beta-\alpha=\pm\pi/2$ (Eq. \eqref{al-be}), and $\beta-\alpha$ never intrudes (Eq. \eqref{abIneq}) into the region $\pm\pi/2$  delineated by dotted lines .}
\label{fig:bma}
\end{figure}

Recall that $f(\rho)<0$ for $\rho\gtrsim6.2$. Equations \eqref{abIneq}, \eqref{al-be} and \eqref{L1Neg} between them provide a reasonable prediction of the structure of the function $\beta-\alpha$ - see Figure \ref{fig:bma}.

\subsection{Inverse of the singular transformation}

Although the inverse transformation of \eqref{ZR} and \eqref{ZI} is singular, it is possible to work from \eqref{Z1R} and \eqref{Z1I} in the limit $\sigma=1/2$. In that limit, taking \eqref{sinA=T} and \eqref{cosA=T} into account, we find (Maple)
\begin{equation}
\mapleinline{inert}{2d}{Zeta1[I] = 2*Zeta[R]/(sin(2*alpha)*f)-(cos(2*alpha)+1)*Zeta1[R]/sin(2*alpha)}{\[\displaystyle {\zeta^{\prime}}_{I}={\frac {2\,\zeta_{{R}}}{\sin \left( 2\,\alpha \right) f}}-{\frac { \left( \cos \left( 2\,\alpha \right) +1 \right) {\it \zeta^{\prime}}_{{R}}}{\sin \left( 2\,\alpha \right) }}\]}
\label{Z1Ia}
\end{equation}
\labeldefs{\begin{flushright}Z1Ia\end{flushright}}\fi
which, when fully converted to polar form, reduces to a tautology after applying \eqref{ZetaRR}. When  \eqref{Z1R} is calculated to first order in $\sigma-1/2$, we find
%taken after Eq. 32 of singular@half.mw%
\begin{equation}
\mapleinline{inert}{2d}{Z1R = (1/2)*fd*ZI/f+ZR/f+((1/2)*sin(alpha)^2*Z2R-(1/4)*sin(2*alpha)*Z2I)*f}{\[\displaystyle {\zeta^{\prime}_R}={\frac {{f^{\prime}}\,{\zeta_I}}{2\,f}}+{\frac {{\zeta_R}}{f}}+ \frac{1}{4}\left( 2\sin^2 \left( \alpha \right)  {\zeta^{\prime\prime}_R}-\sin \left( 2\,\alpha \right) {\zeta^{\prime\prime}_I}
\mbox{} \right) f\]}
\label{Z1RIa}
\end{equation}
\labeldefs{\begin{flushright}Z1RIa\end{flushright}}\fi
and
\begin{equation}
\mapleinline{inert}{2d}{Z1I = -(1/2)*fd*ZR/f+ZI/f+((1/2)*cos(alpha)^2*Z2I-(1/4)*sin(2*alpha)*Z2R)*f}{\[\displaystyle {\zeta^{\prime}_I}=-{\frac {{f^{\prime}}\,{\zeta_R}}{2\,f}}+{\frac {{\zeta_I}}{f}}+ \frac{1}{4}\left( 2\,  \cos^2 \left( \alpha \right)  {\zeta^{\prime\prime}_I}-\sin \left( 2\,\alpha \right) {\zeta^{\prime\prime}_R}
\mbox{} \right) f\]}\,.
\label{Z1RIb}
\end{equation}
\labeldefs{\begin{flushright}Z1RIb\end{flushright}}\fi

An interesting set of results can be obtained from \eqref{Z1RIa} and \eqref{Z1RIb} by  multiplying \eqref{Z1RIa} by $\zeta_I$, \eqref{Z1RIb} by $\zeta_R$, subtracting, converting to polar form and applying \eqref{ZZp}, yielding

\begin{equation}
\mapleinline{inert}{2d}{abs(Zeta[2])/abs(Zeta[1]) = (-fd*cos(alpha-beta)+2*sin(alpha-beta))/(sin(-gamma+alpha)*f)}{\[\displaystyle {\frac { \left| {\zeta^{\prime\prime}} \right| }{ \left| {\zeta^{\prime}} \right| }}={\frac {-{f^{\prime}}\,\cos \left( \alpha-\beta \right) +2\,\sin \left( \alpha-\beta \right) 
\mbox{}}{\sin \left( -\gamma(\rho)+\alpha \right) f}}\]}\,.
\label{Z2ZI}
\end{equation}
\labeldefs{\begin{flushright}Z2ZI\end{flushright}}\fi

A related result can be obtained by making use of \eqref{ZZp}, \eqref{ZpZp} and \eqref{LogZA} to obtain
\begin{equation}
\mapleinline{inert}{2d}{abs(Zeta[2])/abs(Zeta) = (-f^{\prime}+2*tan(alpha-beta))/(sin(-gamma+alpha)*f^2)}{\[\displaystyle {\frac { \left| {\zeta^{\prime\prime}} \right| }{ \left| \zeta \right| }}={\frac {-{f^{\prime}}+2\,\tan \left( \alpha-\beta \right) }{\sin \left( -\gamma(\rho)+\alpha \right) 
\mbox{}{f}^{2}}}\]}\,,
\label{Z2ZIb}
\end{equation}
\labeldefs{\begin{flushright}Z2ZIb\end{flushright}}\fi
where, in both cases, $\gamma$ has been written $\gamma(\rho)$ to distinguish it from Euler's constant (see \eqref{ZDouble}). 
  
\subsection{Complex representations on  the critical line}
By combining the real and imaginary components of $\zeta(1/2+i\rho)$ into a complex representation, we find several interesting forms equivalent to \eqref{ZR} and \eqref{ZI}, those being
\begin{equation}
\mapleinline{inert}{2d}{Zeta(1/2+i*rho) =
1/2*f*(exp(2*i*alpha)*Zeta(1,1/2-i*rho)+Zeta(1,1/2+i*rho));}{%
\[
{\displaystyle \frac{2}{f}}\,\zeta ({\displaystyle  {1}/{2}}  + i\,\rho )=e^{2\,i\,\alpha }\,\zeta^{\prime} (
{\displaystyle {1}/{2}}  - i\,\rho ) + \zeta^{\prime} (
{\displaystyle {1}/{2}}  + i\,\rho )
\]
}
\label{form1}
\end{equation}
\labeldefs{\begin{flushright}form1\end{flushright}}\fi
equivalent to 
\begin{equation}
|\zeta(1/2+i\rho)|=e^{-i\alpha}\zeta(1/2+i\rho)=f\,\Re \left( e^{i\alpha}\zeta^{\prime}(1/2-i\rho)\right)
\label{form2}
\end{equation}
\labeldefs{\begin{flushright}form2\end{flushright}}\fi

and

\begin{equation}
\frac{\zeta(1/2+i\rho)}{|\zeta^{\prime}(1/2+i\rho)|}=f\,e^{i\alpha}\cos(\alpha-\beta)\,
\label{ZabsZp}
\end{equation}
\labeldefs{\begin{flushright}ZabsZp\end{flushright}}\fi

or, equivalently, either
\begin{equation}
\frac{\zeta(1/2+i\rho)}{\zeta^{\prime}(1/2+i\rho)}=fe^{i(\alpha-\beta)}\cos(\alpha-\beta)\,
\label{ZoverZp}
\end{equation}
\labeldefs{\begin{flushright}ZoverZp\end{flushright}}\fi
or \eqref{ZZp}.\newline

Simple expansion of \eqref{form2} into its real and imaginary parts shows that it is equivalent to Eq.(2.12) of I. With respect to the above, the requirement that the right-hand side of \eqref{ZoverZp} must vanish at a zero, is equivalent to Eq.(7.8) of I, provided that $\zeta^{\prime}(1/2+i\rho_0)\neq0$.

\section{At a Zero}
\subsection{Conditions on $\beta(s_0)$}
%see file zeros_on_s.mw
At a zero ($s=s_0$), we require that $\zeta_R(s_0)=\zeta_I(s_0)=\widetilde{\zeta}_R(s_0)=\widetilde{\zeta}_I(s_0)=0$; with these conditions, and solving \eqref{ZAR0} and \eqref{ZAI0} for $|\zeta^{\prime}(s_0|^2$, respectively, we find that, at a zero anywhere in the complex $s$ plane, $|\zeta^{\prime}(s_0|^2$ must satisfy
\begin{equation}
\mapleinline{inert}{2d}{abs(Zeta1)^2 =
2*(Zeta1[I]*Zeta1a[I]-Zeta^{\prime}[R]*\widetilde{zeta}^{\prime}1a[R])*(2*Pi)^sigma/g[2];}{%&
\[
 \left|  \! \,\zeta^{\prime}(s_0)\, \!  \right| ^{2}  =  {\displaystyle \frac {2
\,({\zeta^{\prime}_{I}(s_0)}\,{\mathit{\widetilde{\zeta}^{\prime}}_{I}(s_0)} - {\zeta^{\prime}_{R}(s_0)}\,{
\mathit{\widetilde{\zeta}^{\prime}}_{R}(s_0)})\,(2\,\pi )^{\sigma }}{{g_{2}(s_0)}}} 
\]
\label{Zabs1}
}
\end{equation}
\labeldefs{\begin{flushright}Zabs1\end{flushright}}\fi
and simultaneously 
\begin{equation} 
\mapleinline{inert}{2d}{abs(Zeta1)^2 =
-2*(Zeta1[I]*Zeta1a[R]+Zeta1[R]*Zeta1a[I])*(2*Pi)^sigma/g[1];}{%
\[
 \left|  \! \,\zeta^{\prime}(s_0)\, \!  \right| ^{2}  =  - {\displaystyle 
\frac {2\,({\zeta^{\prime}_{I}(s_0)}\,{\mathit{\widetilde{\zeta}^{\prime}}_{R}(s_0)} + {\zeta^{\prime}_{R}(s_0)}\,
{\mathit{\widetilde{\zeta}^{\prime}}_{I}(s_0)})\,(2\,\pi )^{\sigma }}{{g_{1}(s_0)}}} 
\]
\label{Zabs2}
}
\end{equation}
\labeldefs{\begin{flushright}Zabs2\end{flushright}}\fi

Similarly, from \eqref{ZAR} and \eqref{ZAI} we obtain the simpler conditions
\begin{equation}
\mapleinline{inert}{2d}{T2c := Zeta[1,R] =
-1/8*(2*Pi)^sigma*(Zeta[1,a,I]*g[1]+Zeta[1,a,R]*g[2])/abs(GAMMA(s))^2/
c[0];}{%
\[
{\zeta^{\prime} _{R}(s_0)}= - {\displaystyle \frac {1}{8}
} \,{\displaystyle \frac {(2\,\pi )^{\sigma_0 }\,({{\widetilde{\zeta}}^{\prime}_{ 
\,I}(s_0)}\,{g_{1}(s_0)} + {\widetilde{\zeta}^{\prime}_{\,R}(s_0)}\,{g_{2}(s_0)})}{ \left|  \! \,
\Gamma (s_0)\, \!  \right| ^{2}\,{c_{0}}}} 
\]
}
\label{ZcondR}
\end{equation}
\labeldefs{\begin{flushright}ZcondR\end{flushright}}\fi
and
\begin{equation}
\mapleinline{inert}{2d}{T1c := Zeta[1,I] =
1/8*(2*Pi)^sigma*(Zeta[1,a,I]*g[2]-Zeta[1,a,R]*g[1])/abs(GAMMA(s))^2/c
[0];}{%
\[
{\zeta^{\prime}_{\,I}(s_0)}={\displaystyle \frac {1}{8}} \,
{\displaystyle \frac {(2\,\pi )^{\sigma_0 }\,({\widetilde{\zeta}^{\prime}_{\,I}(s_0)
}\,{g_{2}(s_0)} - {\widetilde{\zeta}^{\prime}_{\,R}(s_0)}\,{g_{1}(s_0)})}{ \left|  \! \,
\Gamma (s_0)\, \!  \right| ^{2}\,{c_{0}}}\,.} 
\]
}
\label{ZcondI}
\end{equation}
\labeldefs{\begin{flushright}ZcondI\end{flushright}}\fi

By setting the two right-hand sides of \eqref{Zabs1} and \eqref{Zabs2} equal, and identifying the appropriate ratios of components, we find a necessary, but not sufficient condition that $\zeta(s_0)=0$, that being
\begin{equation}
\tan(\widetilde{\beta}(s_0))=-\frac{g_2(s_0)\,\tan(\beta(s_0))-g_1(s_0)}{\tan(\beta(s_0))\,g_1(s_0)+g_2(s_0)}.
\label{tanBetaW}
\end{equation}
\labeldefs{\begin{flushright}tanBetaW\end{flushright}}\fi
and its inverse (symmetrical under $\beta(s_0)\leftrightarrow \widetilde{\beta}(s_0)$); i.e.
\begin{equation}
\tan({\beta}(s_0))=-\frac{g_2(s_0)\,\tan(\widetilde{\beta}(s_0))-g_1(s_0)}{\tan(\widetilde{\beta}(s_0))\,g_1(s_0)+g_2(s_0)}.
\label{tanBetaWInv}
\end{equation}
\labeldefs{\begin{flushright}tanBetaWInv\end{flushright}}\fi

The result \eqref{tanBetaW} is not unexpected, since it is equivalent to applying l'H\^{o}pital's rule to \eqref{tans} in the limit $\widetilde{\zeta}_I(s_0)=\widetilde{\zeta}_R(s_0)\rightarrow0$ where the limiting ratio redefines $\tan(\widetilde{\alpha}(s_0))\rightarrow \tan(\widetilde{\beta}(s_0))$, if such points exist; on the critical line $s_0=1/2+i\rho_0$, we have $\widetilde{\beta}= \beta$; thus solving \eqref{tanBetaW} in this case gives

\begin{equation}
\mapleinline{inert}{2d}{tan(beta) = (-g[2]\pm(g[1]^2+g[2]^2)^(1/2))/g[1];}{%
\[
\mathrm{tan}(\beta(\rho_0) )={\displaystyle \frac { - {g_{2}} \pm \sqrt{{g
_{1}}^{2} + {g_{2}}^{2}}}{{g_{1}}}} 
\]
}
\label{tanBeta}
\end{equation}
\labeldefs{\begin{flushright}tanBeta\end{flushright}}\fi

reducing to a known result (Eq. (2.9) of I with $\alpha\rightarrow\beta$) after applying identifications given in the Appendix with $\sigma=1/2$, particularly \eqref{g12^2}. A simpler result is also available, by noting that, as proven in I, setting $\mathfrak{L}(s)=0$ in \eqref{variant} results in a necessary condition that $\zeta(s)=0$. Evaluating the ratio of the real and imaginary parts of $\mathfrak{L}(s)$ and converting to polar form, eventually (Maple) yields a necessary condition for locating $s_0$ corresponding to $\zeta(s_0)=0$, that being
\begin{equation}
\mapleinline{inert}{2d}{tan(beta+beta[1]) =
-(-tan(1/2*Pi*sigma)*tanh(1/2*Pi*rho)+tan(rho[theta]))/(1+tan(1/2*Pi*s
igma)*tan(rho[theta])*tanh(1/2*Pi*rho));}{%
\[
\mathrm{tan}(\beta(s_0)  + {\widetilde{\beta}(s_0)})= - {\displaystyle \frac { - 
\mathrm{tan}({\displaystyle {\pi \,\sigma_{0} }/{2}} )\,\mathrm{
tanh}({\displaystyle  {\pi \,\rho_0 }/{2}} ) + \mathrm{tan}({
\rho _{\theta(s_0) }})}{1 + \mathrm{tan}({\displaystyle  {\pi \,
\sigma_{0} }/{2}} )\,\mathrm{tan}({\rho _{\theta(s_0) }})\,\mathrm{tanh}(
{\displaystyle {\pi \,\rho_0 }/{2}} )}}={\displaystyle \frac{g_1(s_0)}{g_2(s_0)}} \,, 
\]
}
\label{BetaCondition}
\end{equation}
\labeldefs{\begin{flushright}BetaCondition\end{flushright}}\fi
the second equality arising due to \eqref{g12_polar}.\newline

On the critical line, \eqref{BetaCondition} reduces to the following numerical condition for locating $\rho=\rho_0$ corresponding to $\beta=\beta_0$, equivalent to \eqref{al-be}):
\begin{equation}
\mapleinline{inert}{2d}{tan(2*beta(s[0])) =
-(-tanh(1/2*Pi*rho[0])+tan(rho[theta(s[0])]))/(1+tan(rho[theta(s[0])])
*tanh(1/2*Pi*rho[0]));}{%
\[
\mathrm{tan}(2\,\beta_0 )= - {\displaystyle \frac { - 
\mathrm{tanh}({\displaystyle} \,\pi \,{\rho_0/2})
 + \mathrm{tan}({\rho _{\theta_0 }})}{1 + \mathrm{tan}({
\rho _{\theta_0 )}})\,\mathrm{tanh}({\displaystyle} \,\pi \,{\rho_0/2 })}}\,. 
\]
}
\label{Tan2Beta}
\end{equation}
\labeldefs{\begin{flushright}Tan2Beta\end{flushright}}\fi 
Notice that the right hand side of this criterion corresponds to $\tan(2\alpha)$ for all values of $\rho$ because of \eqref{sinA=T} and \eqref{cosA=T}.\newline

Asymptotically ($\rho\rightarrow\infty$), \eqref{BetaCondition} further reduces to 
\begin{equation}
\mapleinline{inert}{2d}{tan(2*beta_0) = tan(1/4*Pi-theta_0+rho_0*ln(2*Pi));}{%
\[
\mathrm{tan}(\beta(s_0)  + {\widetilde{\beta}(s_0)})\sim\mathrm{tan}({\displaystyle \frac {\pi\sigma_0 }{
2}}  - \rho_{\theta(s_{0})} ))\,.
\]
}
\label{Tan2Beta_Asy_Gen}
\end{equation}
\labeldefs{\begin{flushright}$Tan2Beta{_{Asy}}_{Gen}$\end{flushright}}\fi 
and \eqref{Tan2Beta} becomes
\begin{equation}
\mapleinline{inert}{2d}{tan(2*beta_0) = tan(1/4*Pi-theta_0+rho_0*ln(2*Pi));}{%
\[
\mathrm{tan}(2\,\beta_0 )\sim\mathrm{tan}({\displaystyle \frac {\pi }{
4}}  - \rho_{\theta_0}  )\,.
\]
}
\label{Tan2Beta_Asy}
\end{equation}
\labeldefs{\begin{flushright}$Tan2Beta_Asy$\end{flushright}}\fi

\section{Summary and a Comment}
In the previous sections, new relationships have been developed among the real and imaginary components of $\zeta(s)$ on both sides of the critical strip. Notably, an analytic expression was obtained for the sum $\alpha(s)+\widetilde{\alpha}(s)$ as well as a simplified form of the transformation between the real and imaginary components of $\zeta$ and $\zeta^{\prime}$ on the critical strip.  In addition, simplified derivations of previous results were both presented and generalized from the critical strip to the whole complex plane.\newline 

The results relating various functions on both sides of the critical strip are significant because of interest in several theorems in the literature, generally based on an analysis of \eqref{Zmags}. Spira \cite{Spira1965}, Saidak and Zvengrowski \cite{Sai&Zven}, Nazardonyavi and Yakubovich \cite{Naz&Yaku} show for $0<\sigma<1/2$, that 
\begin{equation}
|\zeta(s)|\geq |\widetilde{\zeta}(s)|,
\label{inequal}
\end{equation}
\labeldefs{\begin{flushright}$inequal$\end{flushright}}\fi
``with equality only if $|\zeta(s)|=0$". All of these claim that if the ``$\geq$" operator could be replaced by the ``$>$" operator, RH would be proven. %(\cite{[Spira1965,Sai&Zven,Naz&Yaku,ZetaNotes})% 
These claims are incorrect.\newline

Consider the function (Milgram, see Eq. (A.1)) \cite {Milgram}
% from file counterexample.mw%
\begin{equation}
\mapleinline{inert}{2d}{Zeta[c] := Zeta(s)*sin(Pi*(s-s0))*sin(Pi*(s+s0))*sin(Pi*(s-conjugate(s0)))*sin(Pi*(s+conjugate(s0)))/(cosh(Pi*rho0)^2-cos(Pi*sigma0)^2)^2}{\[\displaystyle \zeta_{c}(s)\, \equiv\,{\frac {\zeta \left( s \right) \sin \left( \pi\, \left( s-{ s_0} \right)  \right) \sin \left( \pi\, \left( s+{ s_0} \right)  \right) \sin \left( \pi\, \left( s-\overline{{ s_0}} \right)  \right) 
\mbox{}\sin \left( \pi\, \left( s+\overline{{ s_0}} \right)  \right) }{ \left(  \cosh^2 \left( \pi\,\rho_0   \right) -  \cos^2 \left( \pi\,\sigma_0 \right)  \right)  ^{2}}}\]}
\label{zeta_C}
\end{equation}
\labeldefs{\begin{flushright}$zeta_C$\end{flushright}}\fi
where $s_0\equiv \sigma_0+i\rho_0$ with $\rho_0>1$ locates a (complex, arbitrary) zero and $\overline{{ s_0}}$ denotes complex conjugation. This function has the interesting properties that it satisfies the functional equation \eqref{basic} because
\begin{equation}
\frac{\zeta_c(s)}{\zeta_c(1-s)}=\frac{\zeta(s)}{\zeta(1-s)}\,,
\label{Zcrat}
\end{equation}
\labeldefs{\begin{flushright}Zcrat\end{flushright}}\fi
and, in common with $\zeta(s)$, is both self-conjugate (see \eqref{cmplxC}), and possesses a pole with residue unity at $s=1$. As noted, $\zeta_c$ possesses complex zeros at $s=s_0,s=1-s_0$ and conjugate points, but cannot be confused with $\zeta(s)$ because it also possesses zeros at $s=s_0\pm\,k$, where it is well-known, $\zeta(s)$ does not. In fact, any function (but, see Titchmarsh and Heath-Brown, \cite[Section 2.13]{Titch2}) of the form
\begin{equation}
\Upsilon(s)=w(s)\zeta(s)
\label{Updef}
\end{equation}
\labeldefs{\begin{flushright}Updef\end{flushright}}\fi
will satisfy the functional equation provided that $w(s)$ is self-conjugate and satisfies
\begin{equation}
w(s)=w(1-s).
\label{ws}
\end{equation}
\labeldefs{\begin{flushright}ws\end{flushright}}\fi

% from file absZratio.mw%
It is a simple matter to recognize that at a complex zero ($s=s_0$) of order $n$, the expression \eqref{Zmags} reduces to the well-defined limit
\begin{equation}
\lim_{s \to s_0}\,\displaystyle {\frac { \displaystyle | \zeta(s) |^{2}}{  | \widetilde{\zeta}(s)|^{2}}}
=\displaystyle {\frac { (\displaystyle | \zeta(s_0) |^{2})^{(n)}}{(  | \widetilde{\zeta}(s_0)|^{2})^{(n)}}}
=\Phi(s_0)\,
\label{Zmags0}
\end{equation}
\labeldefs{\begin{flushright}Zmags0\end{flushright}}\fi
in terms of derivatives of $n^{th}$ order, by l'H\^opital's rule applied to an analytic function, and thus none of the cited allusions to a possible ``proof" of RH are valid. A simple study of \eqref{zeta_C} exemplifies this reasoning. For $0\leq\sigma<1/2$, evaluating the simple limit, yields 

\begin{equation}
 \lim_{s\to s_0} \displaystyle {\frac { \displaystyle | \zeta_c(s) |^{2}}{  | \widetilde{\zeta_c}(s)|^{2}}} = \displaystyle {\frac { \displaystyle | \zeta(s_0) |^{2}}{  | \widetilde{\zeta}(s_0)|^{2}}}>1,
 \label{ratios}
\end{equation}
\labeldefs{\begin{flushright}ratios\end{flushright}}\fi

where the inequality (lack of equality) is implied by any of the cited proofs, because $s_0$ is arbitrary and, without loss of generality, we can specify that $\zeta(s_0)\ne0$. To reiterate, the function $\zeta_c(s)$ demonstrates that it is possible for any function that satisfies the functional equation and is self-conjugate to possess a complex zero, satisfy the stronger form of \eqref{inequal} and violate RH. 
\newline

\begin{figure}[h] % from file "counterexample2.mw"%
\centering
\includegraphics[width=.41\textwidth]{{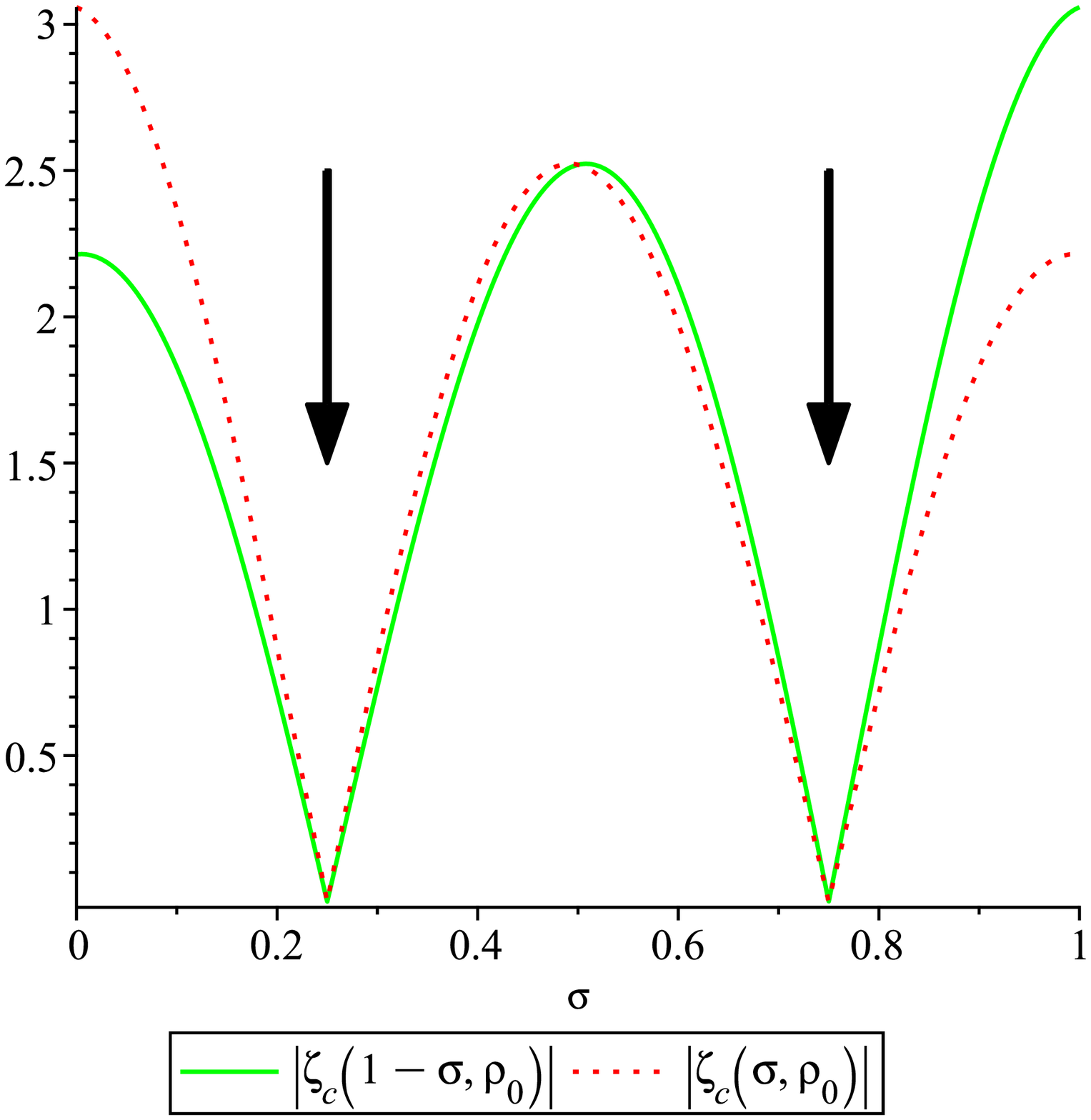}}
\includegraphics[width=.4\textwidth]{{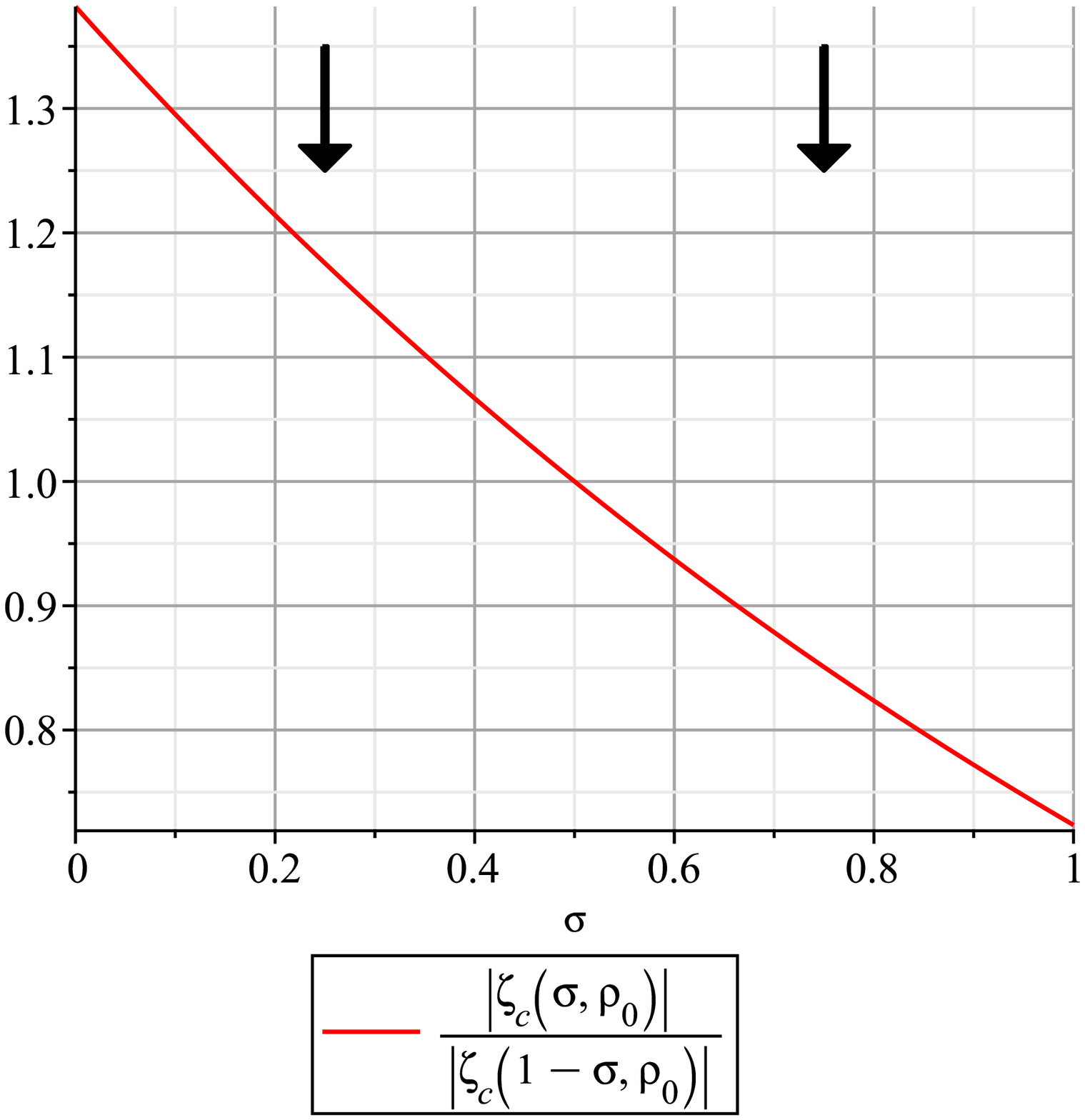}}
\caption{Plot of $|\zeta_c(\sigma,\rho_0)|$ and $|\zeta_c(1-\sigma,\rho_0)|$ {\bf (left)} and $|\zeta_c(\sigma,\rho_0)|/|\zeta_c(1-\sigma,\rho_0)|$ {\bf (right)} as a function of $\sigma$ using $s_0=3/4+i\rho_0$ with $\rho_0=12$, demonstrating the inequality \eqref{ratios} at $s_0$ in the limiting case $s=s_0$. The arrows point to the location of zeros.} 
\label{fig:countfig}
\end{figure}

Further insight on this subject can be obtained by straightforward evaluation of the derivative of $\Phi(s)$ with respect to $\sigma$ (Maple), giving (and providing a simplified derivation for results obtained by Nazardonyavi and Yakubovich \cite{Naz&Yaku})
\begin{equation}
\mapleinline{inert}{2d}{Diff(abs(Zeta(sigma+I*rho))^2/abs(Zeta(1-sigma+I*rho))^2, sigma) = (1/2)*(2*Pi)^(2*sigma)*(Pi*sin(Pi*sigma)-Psi[2]*(cos(Pi*sigma)+cosh(Pi*rho)))/((cos(Pi*sigma)+cosh(Pi*rho))^2*abs(GAMMA(sigma+I*rho))^2)}{\[\displaystyle {\frac {\partial }{\partial \sigma}} \left( {\frac { \left(  \left| \zeta \left( \sigma+i\rho \right)  \right|  \right) ^{2}}{ \left(  \left| \zeta \left( 1-\sigma+i\rho \right)  \right|  \right) ^{2}}} \right) ={\frac { \left( 2\,\pi \right) ^{2\,\sigma}
\mbox{} \left( \pi\,\sin \left( \pi\,\sigma \right) -\Psi_{{2}} \left( \cos \left( \pi\,\sigma \right) +\cosh \left( \pi\,\rho \right)  \right)  \right) }{ 2\,\left( \cos \left( \pi\,\sigma \right) +\cosh \left( \pi\,\rho \right)  \right) ^{2}
\mbox{} \left(  \left| \Gamma \left( \sigma+i\rho \right)  \right|  \right) ^{2}}}\]}
\label{ZmagsDiff}
\end{equation}
\labeldefs{\begin{flushright}ZmagsDiff\end{flushright}}\fi

from which we conclude that $\displaystyle {\frac { \displaystyle | \zeta(s) |^{2}}{  | \widetilde{\zeta}(s)|^{2}}}>1$ because the right-hand side of \eqref{ZmagsDiff} is obviously always negative and monotonic for $\rho\gtrsim 10$ and all $\sigma$ (see \eqref{Psi2} and \eqref{PsiRAsy}). Therefore $\Phi(s)>1$ when $0\leq\sigma<1/2$ independent of the possibility that $\zeta(s_0)=0$. This further demonstrates the invalidity of the comments cited, including some made by myself  in \cite{ZetaNotes}. See also Nazardonyavi and Yakubovich \cite[Proposition (2)]{Naz&Yaku}. For further clarity, see Figure \ref{fig:countfig}.

\section{Appendix- notation}

The following summarizes the notation used throughout. Each of the symbols retains functional dependence according to how it is referenced in the main text. Any symbol lacking specific functional dependence is assumed to be only a function of $\rho$. If some entity is referenced only as a function of $\rho$, the implication is that any appearance of $\sigma$ in its structure corresponds to $\sigma=1/2$. Subscripts $R$ and $I$ respectively refer to the real and imaginary components of the associated quantity. In all cases $\rho\geq0$, $k$ is an integer and $\psi$ refers to the digamma function. All derivatives specified by the ``prime" symbol ($^{\prime}$) are taken with respect to $\rho$ unless specified otherwise. In polar notation, I use

\begin{eqnarray}
\zeta(s)=e^{i\alpha(s)}|\zeta(s)|
\label{Zeta-polar}\\
\widetilde{\zeta}(s)=e^{i\widetilde{\alpha}(s)}|\widetilde{\zeta}(s)|
\label{Zeta-polar-s}
\end{eqnarray}
\labeldefs{\begin{flushright}Zeta-polar,Zeta-polar-s\end{flushright}}\fi
\begin{eqnarray}
\Gamma(s)=e^{i\theta(s)}|\Gamma(s)| \\
\zeta^{\prime}(s)=e^{i\beta(s)}|\zeta^{\prime}(s)|\\
\widetilde{\zeta}^{\prime}(s)=e^{i\widetilde{\beta}(s)}|\widetilde{\zeta}^{\prime}(s)|\\  
\zeta^{\prime\prime}(s)=e^{i\gamma(s)}|\zeta^{\prime\prime}(s)| \label{ZDouble} \\
\widetilde{\zeta}^{\prime\prime}(s)=e^{i\widetilde{\gamma}(s)}|\widetilde{\zeta}^{\prime\prime}(s)|  
\label{zetaExp2}
\end{eqnarray}

For example, specific to the ``critical line" $s=1/2+i\rho$, a specialized form might be written
\begin{equation}
\zeta(1/2+i\rho)=e^{i\alpha}|\zeta|=e^{i\alpha}\sqrt{\zeta_{R}^2+\zeta_{I}^2}
\label{ZetaExp}
\end{equation}
\labeldefs{\begin{flushright}ZetaExp\end{flushright}}\fi

with symbols

\begin{eqnarray}
\rho_{\pi}=\rho\log(2\pi)
\label{rhoPi}
\\
\theta={arg}(\Gamma(1/2+i\rho))+k \pi
\label{argGamma}
\\
\alpha={arg}(\zeta(1/2+i\rho))+k \pi
\label{argZeta}
\\
\beta={arg}(\zeta^{\prime}(1/2+i\rho))+k \pi
\label{argZeta1}
\\
\gamma={arg}(\zeta^{\prime\prime}(1/2+i\rho))+ k\pi
\label{argZeta2}\,.
\end{eqnarray}
\labeldefs{\begin{flushright}rhoPi,rhoTheta,argGamma,argZeta\end{flushright}}\fi
The symbols $(\theta,\alpha,\beta,\gamma)$ are continuous functions, whereas the $arg$ operator denotes the corresponding discontinuous function limited to $(-\pi,\pi)$, the two being separated by a term equal to $k\pi$, $k=0,\pm 1,\pm 2,\dots$.\newline

In general
\begin{eqnarray}
&\rho_{\theta(s)}=\rho_{\pi}-\theta(s)
\label{rhoTheta}
\\
&\theta(s)={arg}(\Gamma(\sigma+i\rho))+ k\pi=\Im({\rm Log\Gamma}(s))
\label{argGammaS}
\\
&\alpha(s)={arg}(\zeta(\sigma+i\rho))+k\pi
\label{argZetaS}
\\
&{\mathrm {and}}\\
&\widetilde{\alpha}(s)={arg}(\zeta(1-\sigma+i\rho)) +k\pi\,.
\end{eqnarray}
\labeldefs{\begin{flushright}argGammaS,argZetaS\end{flushright}}\fi

Specialized symbols are
\begin{eqnarray}
&c_{{0}}=1/2\,\cos \left( \pi \,\sigma \right) +1/2\,\cosh \left( \pi 
\,\rho \right) \label{c_0} \\ 
&g_{{1}}(s)=4\,\Gamma _{{I}}(s)S_{{2}}(s)+4\,\Gamma _{{R}}(s)S_{{1}}(s)  \label{g_1} \\
&g_{{2}}(s)=4\,\Gamma _{{I}}(s)S_{{1}}(s)-4\,\Gamma _{{R}}(s)S_{{2}}(s) \label{g_2} \\
&S_{{1}}(s)=\sin \left( \rho_{{\pi }} \right) \cos \left( 1/2\,\pi \,
\sigma \right) \cosh \left( 1/2\,\pi \,\rho \right) +\cos \left( \rho_
{{\pi }} \right) \sin \left( 1/2\,\pi \,\sigma \right) \sinh \left( 1/
2\,\pi \,\rho \right) \label{s_1} \\ 
&S_{{2}}(s)=\sin \left( \rho_{{\pi }} \right) \sin \left( 1/2\,\pi \,
\sigma \right) \sinh \left( 1/2\,\pi \,\rho \right) -\cos \left( \rho_
{{\pi }} \right) \cos \left( 1/2\,\pi \,\sigma \right) \cosh \left( 1/
2\,\pi \,\rho \right) \label{s_2}
\end{eqnarray}
\labeldefs{\begin{flushright}$c_0,c_1,g_1,g_2,s_1,s_2,psi_2$\end{flushright}}\fi
\begin{equation}
\mapleinline{inert}{2d}{f(s) = ln(2*Pi)-psi(s)+1/2*Pi*tan(1/2*Pi*s);}{%
\[
\mathit{f}(s)=\mathrm{ln}(2\,\pi ) - \psi (s) + {\displaystyle 
\frac {1}{2}} \,\pi \,\mathrm{tan}({\displaystyle \frac {\pi \,s
}{2}} )
\]
}
\label{fs}
\end{equation}
\labeldefs{\begin{flushright}fs\end{flushright}}\fi
\begin{equation}
\mapleinline{inert}{2d}{f = 4*cosh(Pi*rho)/(2*ln(2*Pi)*cosh(Pi*rho)-2*Re(Psi(1/2+I*rho))*cosh(Pi*rho)+Pi)}{\[\displaystyle f={\frac {4\,\cosh \left( \pi\,\rho \right) }{2\,\ln  \left( 2\,\pi \right) \cosh \left( \pi\,\rho \right) -2\,\Re \left( \psi \left( 1/2+i\rho \right)  \right) \cosh \left( \pi\,\rho \right) 
\mbox{}+\pi}}\]}
\label{f}
\end{equation}
\labeldefs{\begin{flushright}f\end{flushright}}\fi
\begin{eqnarray}
g_1(s)^2+g_2(s)^2=16\,|\Gamma(s)|^2\,c_0 \\
g_1^2+g_2^2=8\pi
\label{g12^2}
\end{eqnarray}
\labeldefs{\begin{flushright}$g12^2$\end{flushright}}\fi
In polar form
\begin{eqnarray} \nonumber
\mapleinline{inert}{2d}{g[1] :=
(4*sin(1/2*Pi*sigma)*cos(theta-rho*ln(2*Pi))*sinh(1/2*Pi*rho)-4*cos(1/
2*Pi*sigma)*sin(theta-rho*ln(2*Pi))*cosh(1/2*Pi*rho))*abs(GAMMA);}{%
\maplemultiline{
{g_{1}(s)}=4\, \left|  \! \,\Gamma(s) \, \!  \right| \left(\mathrm{sin}({\displaystyle \frac {\pi \,\sigma }{
2}} )\,\mathrm{cos}(\rho_{\theta(s)})\,
\mathrm{sinh}({\displaystyle \frac {\pi \,\rho }{2}} ) - 
\mathrm{cos}({\displaystyle \frac {\pi \,\sigma }{2}} )\,\mathrm{
sin}(\rho_{\theta(s)})\,\mathrm{cosh}(
{\displaystyle \frac {\pi \,\rho }{2}} )\right) 
 }
}\\ 
\label{g12_polar}
\mapleinline{inert}{2d}{g[1] = (2*(sin((1/2)*Pi*sigma-rho[theta(s)])*exp((1/2)*Pi*rho)-sin((1/2)*Pi*sigma+rho[theta(s)])*exp(-(1/2)*Pi*rho)))*abs(GAMMA(s))}{\[\displaystyle =2\, \left| \Gamma \left( s \right)  \right| \, \left( \sin \left( \pi\,\sigma/2-\rho_{{\theta \left( s \right) }} \right) {{\rm e}^{1/2\,\pi\,\rho}}
\mbox{}-\sin \left( \pi\,\sigma/2+\rho_{{\theta \left( s \right) }} \right) {{\rm e}^{-1/2\,\pi\,\rho}} \right) \]}\\
\mapleinline{inert}{2d}{g[2] :=
(4*sin(1/2*Pi*sigma)*sin(theta-rho*ln(2*Pi))*sinh(1/2*Pi*rho)+4*cos(1/
2*Pi*sigma)*cos(theta-rho*ln(2*Pi))*cosh(1/2*Pi*rho))*abs(GAMMA);}{%
\maplemultiline{
{g_{2}(s)} = 4\, \left|  \! \,\Gamma(s) \, \!  \right|\left(\mathrm{sin}({\displaystyle \frac {\pi \,\sigma }{
2}} )\,\mathrm{sin}(\rho_{\theta(s)})\,
\mathrm{sinh}({\displaystyle \frac {\pi \,\rho }{2}} ) + 
\mathrm{cos}({\displaystyle \frac {\pi \,\sigma }{2}} )\,\mathrm{
cos}(\rho_{\theta(s)})\,\mathrm{cosh}(
{\displaystyle \frac {\pi \,\rho }{2}} )\right)
  }
}\\
\mapleinline{inert}{2d}{g[2] = 2*abs(GAMMA)*(cos((1/2)*Pi*sigma-rho[theta(s)])*exp((1/2)*Pi*rho)+cos((1/2)*Pi*sigma+rho[theta(s)])*exp(-(1/2)*Pi*rho))}{\[\displaystyle =2\, \left| \Gamma(s) \right|  \left( \cos \left( \pi\,\sigma/2-\rho_{{\theta \left( s \right) }} \right) {{\rm e}^{1/2\,\pi\,\rho}}
\mbox{}+\cos \left( \pi\,\sigma/2+\rho_{{\theta \left( s \right) }} \right) {{\rm e}^{-1/2\,\pi\,\rho}} \right) \]\,,}
\end{eqnarray}
\labeldefs{\begin{flushright}$g12_{polar}$\end{flushright}}\fi

along with the definitions
% see file Check_Eqns1.mw
\begin{eqnarray}
\zeta_p(s)\equiv\zeta_I(s)\,g_1(s)+\zeta_R(s)\,g_2(s)\\
\zeta_m(s)\equiv\zeta_I(s)\,g_2(s)-\zeta_R(s)\,g_1(s)
\label{Zpm}
\end{eqnarray}
\labeldefs{\begin{flushright}Zpm\end{flushright}}\fi 
\begin{equation}
a\equiv \frac{\sqrt{2}\,g_1}{8\,\sqrt{\pi}}
\label{a-def}
\end{equation}
\labeldefs{\begin{flushright}$a-def$\end{flushright}}\fi 
\begin{equation}
b\equiv \frac{\sqrt{2}\,g_2}{8\,\sqrt{\pi}}
\label{b-def}
\end{equation}
\labeldefs{\begin{flushright}$b-def$\end{flushright}}\fi
\begin{equation}
h_1(s)\equiv - {g_{1}(s)}\,{p_{1}(s)} + {g_{2}(s)}\,{p_{2}(s)}
\label{h1s}
\end{equation}
\labeldefs{\begin{flushright}$h1s$\end{flushright}}\fi
%for h1 and h2 see file "singular@half.mw"%
\begin{equation}
h_2(s)\equiv \:\:\:{g_{1}(s)}\,{p_{2}(s)} + {g_{2}(s)}\,{p_{1}(s)}
\label{h2s}
\end{equation}
\labeldefs{\begin{flushright}h2s\end{flushright}}\fi
\begin{equation}
\mapleinline{inert}{2d}{h__3(s) = -g[1](s)*p[2](s)+g[2](s)*p[1](s)}{\[\displaystyle h_{3} \left( s \right) \equiv-g_{1} \left( s \right) p_{{2}} \left( s \right) +g_{{2}} \left( s \right) p_{{1}} \left( s \right) \]}
\label{h3s}
\end{equation}
\labeldefs{\begin{flushright}h3s\end{flushright}}\fi
\begin{equation}
\mapleinline{inert}{2d}{h__4(s) = g[1](s)*p[1](s)+g[2](s)*p[2](s)}{\[\displaystyle h_{4} \left( s \right) \equiv\:\:\:g_{1} \left(s\right) p_{{1}} \left( s \right) +g_{{2}} \left( s \right) p_{{2}} \left( s \right) \]}
\label{h4s}
\end{equation}
\labeldefs{\begin{flushright}h4s\end{flushright}}\fi

and
\begin{equation}
h_1=-\frac{16\sqrt{2\,\pi}\,\cosh(\pi\rho)\cos(2\alpha)}{f}
\label{h1}
\end{equation}
\labeldefs{\begin{flushright}$h1$\end{flushright}}\fi

\begin{equation}
h_2=-\frac{16\sqrt{2\,\pi}\,\cosh(\pi\rho)\sin(2\alpha)}{f}\,.
\label{h2}
\end{equation}
\labeldefs{\begin{flushright}h2\end{flushright}}\fi
Both of the above two results were obtained using \eqref{sinA=T} and \eqref{cosA=T}. Additionally,

%the following were checked at the bottom of file EqGen_5.mw%

\begin{equation}
h_1(s)^2+h_2(s)^2=(g_1(s)^2+g_2(s)^2)\,(p_1(s)^2+p_2(s)^2)\,.
\label{hSq}
\end{equation}
\labeldefs{\begin{flushright}hSq\end{flushright}}\fi
Furthermore, it is convenient to introduce
\begin{eqnarray}
\mapleinline{inert}{2d}{p[1] = 4*(cos(Pi*sigma)+cosh(Pi*rho))*Psi[1]-2*Pi*sinh(Pi*rho);}{%
\[
{p_{1}(s)}=8\,c_0\,{\Psi _{1}} - 2\,\pi \,\mathrm{sinh}(\pi \,\rho )
\]
}
\label{p1}\\
\mapleinline{inert}{2d}{p[2] = (2*cos(Pi*sigma)+2*cosh(Pi*rho))*Psi[2]-2*Pi*sin(Pi*sigma);}{%
\[
{p_{2}(s)}=4\,c_0\,{\Psi _{2}} - 2\,\pi \,\mathrm{sin}(\pi \,\sigma )
\]
}
\label{p2}\\
\iffalse
\mapleinline{inert}{2d}{-p[1]^2+p[2]^2 = -4*(Pi*sinh(Pi*rho)-4*Psi[I]*c[0])^2+4*(Pi*sin(Pi*sigma)-2*Psi[2]*c[0])^2}{\[\displaystyle {p_{{2}}}^{2}-{p_{{1}}}^{2}=-4\, \left( \pi\,\sinh \left( \pi\,\rho \right) -4\,\Psi_{{1}}c_{{0}} \right) ^{2}+4\, \left( \pi\,\sin \left( \pi\,\sigma \right) -2\,\Psi_{{2}}c_{{0}} \right) ^{2}\]}
\\
\fi
\\\pi_{\sigma}=(2\pi)^\sigma\\
\mapleinline{inert}{2d}{Psi[1] = psi_I(s);}{%
\[
{\Psi _{1}}=\psi_I(s)
\]
\label{Psi1}
}\\
\mapleinline{inert}{2d}{Psi[2] = -2*ln(2*Pi)+2*Re(Psi(s));}{%
\[
{\Psi _{2}}= - 2\,\mathrm{ln}(2\,\pi ) + 2\,\psi_R (s)
\]
\label{Psi2}
}\\
\mapleinline{inert}{2d}{q[1] := -Psi[2]+1/2*Pi*sin(Pi*sigma)/c[0];}{%
\[
{q_{1}} =  - {\Psi _{2}} + {\displaystyle \frac {\pi }{2}} \,
{\displaystyle \frac {\mathrm{sin}(\pi \,\sigma )}{{c_{0}}}
} 
\]
\label{q1}
}\\
\mapleinline{inert}{2d}{q[2] := 2*PsI-1/2*sinh(Pi*rho)*Pi/c[0];}{%
\[
{q_{2}} = 2\,\Psi _{1} - {\displaystyle \frac {\pi}{2}} \,
{\displaystyle \frac {\mathrm{sinh}(\pi \,\rho ) }{{c_{0}}}\,.
} 
\]
}
\label{q2}
\end{eqnarray}

Asymptotically ($\rho\rightarrow\infty, 0<\sigma<1 $) we have \cite[NIST, Eq. 5.11.2]{NIST}
\begin{eqnarray}
\mapleinline{inert}{2d}{PsRe := ln(rho)+1/4*(2*sigma^2-4*sigma+1)/rho^2;}{%
\[
\Re(\psi(\sigma+i\rho)) \approx  \mathrm{ln}(\rho ) + {\displaystyle \frac {2\,
\sigma ^{2} - 4\,\sigma  + 1}{4\,\rho ^{2}}} 
\]
\label{PsiRAsy}
}\\
\mapleinline{inert}{2d}{PsIm :=
1/2*Pi+(-sigma+1)/rho+1/6*sigma*(2*sigma^2-6*sigma+3)/rho^3;}{%
\[
\Im(\psi(\sigma+i\rho)) \approx  {\displaystyle \frac {\pi }{2}}  + 
{\displaystyle \frac {1 - \sigma }{\rho }}  + {\displaystyle 
\frac {\sigma \,(2\,\sigma ^{2} - 6\,\sigma  + 3)}{6\,\rho ^{3}}
} 
\]
}
\label{PsiIAsy}
\end{eqnarray}
\labeldefs{\begin{flushright}PsiRAsy,PsiIAsy\end{flushright}}\fi

whence, in the same asymptotic limit

% from file "Asympt_on_Crit.mw"
\begin{eqnarray}
p_{1} \approx  \left( - {\displaystyle \frac {2\,(\sigma  - 1)}{
\rho }}  + {\displaystyle \frac {\sigma \,(2\,\sigma ^{2} - 6\,
\sigma  + 3)}{3\,\rho ^{3}}} \right) \,e^{\pi \,\rho } + 2\,\pi \,
\mathrm{cos}(\pi \,\sigma )\\ - {\displaystyle \frac {4\,(\sigma 
 - 1)\,\mathrm{cos}(\pi \,\sigma )}{\rho }}  
  +{\displaystyle \frac {2\,\sigma\,\,\mathrm{cos}(\pi \,
\sigma )}{3}} \,{\displaystyle \frac {
 \,(2\,\sigma ^{2} - 6\,\sigma  + 3)}{\rho ^{3}}}  
 \label{p1Asy}
\\
p_{2}\approx \left(  \! \mathrm{ln}({\displaystyle \frac {
\rho ^{2}}{4\,\pi ^{2}}} ) + {\displaystyle \frac {\sigma ^{2} - 
2\,\sigma  + {\displaystyle \frac {1}{2}} }{\rho ^{2}}}  \! 
 \right) \,e^{\pi \,\rho } + \mathrm{cos}(\pi \,\sigma )\,
\mathrm{ln}({\displaystyle \frac {\rho ^{4}}{16\,\pi ^{4}}} )\\
- 2\,\pi \,\mathrm{sin}(\pi \,\sigma ) 
+ {\displaystyle \frac {(2\,\sigma ^{2} - 4\,\sigma  + 1)
\,\mathrm{cos}(\pi \,\sigma )}{\rho ^{2}}} \,.
\label{p2Asy}
\end{eqnarray}
\labeldefs{\begin{flushright}p1Asy,p2Asy\end{flushright}}\fi

Finally, the formal, and straightforwardly obtained (Maple), expressions for the real and imaginary parts of $\widetilde{\zeta}(s)$ as obtained from \eqref{variant} are:

\begin{equation}
\mapleinline{inert}{2d}{ZaR =
((-4*PsI*Z1R+Z1I*Psi[2])*Zeta[g]+8*(-Z1I*Z1aI+Z1R*Z1aR)*(2*Pi)^sigma+4
*g[2]*abs(Zeta1)^2)*CC+(sinh(Pi*rho)*Z1R-sin(Pi*sigma)*Z1I)*Pi*Zeta[g]
;}{%
\maplemultiline{
\mathit{\widetilde{\zeta}_R(s)} = \left( [( - 4\,\mathit{\psi_I(s)}\,\mathit{\zeta^{\prime}_R(s)} + \mathit{\zeta^{\prime}_I}(s)\,{
\Psi _{2}(s)})\,{\zeta _{m}(s)} + 8\,( - \mathit{\zeta^{\prime}_I}(s)\,\mathit{\widetilde{\zeta}^{\prime}_I}(s) + 
\mathit{\zeta^{\prime}_R}(s)\,\mathit{\widetilde{\zeta}^{\prime}_R}(s))\,(2\,\pi )^{\sigma } \right. \\ \left. 
 \hspace{1cm}+\; 4\,{g_{2}(s)}\,
 \left|  \! \,\zeta^{\prime}(s)\, \!  \right| ^{2}]\,\mathit{c_0} 
\mbox{} + (\mathrm{sinh}(\pi \,\rho )\,\mathit{\zeta^{\prime}_R(s)} - \mathrm{sin
}(\pi \,\sigma )\,\mathit{\zeta^{\prime}_I(s)})\,\pi \,{\zeta _{m}(s)} \right)/\zeta_d(s) }
}
\label{ZAR0}
%\end{split}
\end{equation}
\labeldefs{\begin{flushright}ZAR0\end{flushright}}\fi

\begin{equation}
\mapleinline{inert}{2d}{ZaI =
((4*Im(psi(s))*Zeta(1,s)[R]-2*Zeta(1,s)[I]*psi[2])*Zeta[p]+(-8*Zeta(1,
s)[I]*eta(1,s)[R]-8*Zeta(1,s)[R]*eta(1,s)[I])*(2*Pi)^sigma-4*g[1]*abs(
Zeta1)^2)*c[0]+(-sinh(Pi*rho)*Zeta(1,s)[R]+sin(Pi*sigma)*Zeta(1,s)[I])
*Pi*Zeta[p];}{%
\maplemultiline{
\mathit{\widetilde{\zeta}_I(s)} =\,-\,\left( \left[ (4\,\psi_I (s)\,{\zeta^{\prime}_{R}(s)} - 2\,{
\zeta^{\prime}_{I}(s)}\,{\Psi _{2}(s)})\,{\zeta _{p}(s)} 
\mbox{} - 8\,( {\zeta^{\prime}_{I}(s)}\,{\widetilde{\zeta^{\prime}}_{R}(s)} + {
\zeta^{\prime}_R(s)}\,{\widetilde{\zeta^{\prime}}_I(s)})\,(2\,\pi )^{\sigma } \right. \right.  \\ \left.\left.  
\hspace{1.2cm} -\,4\,{g_{1}}\, \left|  \! \,\zeta^{\prime}(s) \, \!  \right| ^{2}\,\right]\,{c_{0}} 
\mbox{} - \left[ \mathrm{sinh}(\pi \,\rho )\,{\zeta^{\prime}_R (s)} - 
\mathrm{sin}(\pi \,\sigma )\,{\zeta^{\prime}_I(s)}\right] \,\pi \,{\zeta 
_{p}(s)}\right)/\zeta_d(s) }
}
\label{ZAI0}
\end{equation}
\labeldefs{\begin{flushright}ZAI0\end{flushright}}\fi
where

\begin{equation}
\mapleinline{inert}{2d}{zeta[d] :=
((-8*Im(Psi(s))*c+2*Pi*sinh(Pi*rho))*Z1I+(-4*c*Psi[2]+2*Pi*sin(Pi*sigm
a))*Z1R)*(2*Pi)^sigma;}{%
\[
{\zeta _{d}(s)} \equiv (( - 8\,\mathit{\psi_I (s)}\,c_0 + 2\,\pi \,\mathrm{sinh
}(\pi \,\rho ))\,\mathit{\zeta^{\prime}_{I}(s)} + ( - 4\,c_0\,{\Psi _{2}(s)} + 2\,\pi \,
\mathrm{sin}(\pi \,\sigma ))\,\mathit{\zeta^{\prime}_{R}(s)})\,(2\,\pi )^{\sigma } \,.
\]
}
\label{Zeta_d}
\end{equation}
\labeldefs{\begin{flushright}$Zeta_d$\end{flushright}}\fi 

\bibliographystyle{unsrt}
%\bibliography{c://Physics//Glasser//Biblio}
\bibliography{.//biblio}

\begin{thebibliography}{10}

\bibitem{Milgram_Exploring}
Michael Milgram.
\newblock Exploring {R}iemann's functional equation.
\newblock {\em Cogent Mathematics}, 3(1):1179246, 2016.
\newblock http://dx.doi.org/10.1080/23311835.2016.1179246.

\bibitem{Spira1965}
Robert Spira.
\newblock An inequality for the {R}iemann {Z}eta function.
\newblock {\em Duke Math. J.}, 32(2):247--250, 06 1965.

\bibitem{NIST}
F.~W.~J. Olver, D.~W. Lozier, R.~F. Boisvert, and C.~W. Clark, editors.
\newblock {\em {NIST Handbook of Mathematical Functions}}.
\newblock Cambridge University Press, New York, NY, 2010.
\newblock Print companion to \cite{NIST:DLMF}.

\bibitem{Spira1973}
Robert Spira.
\newblock Zeros of $\zeta^{\prime} (s)$ and the {R}iemann {H}ypothesis.
\newblock {\em Illinois J. Math.}, 17(1):147--152, 03 1973.

\bibitem{Maple}
Maplesoft, a division of Waterloo Maple Inc.
\newblock {\em Maple.}

\bibitem{Srin&Zven}
G.K. Srinivasan and P.~Zvengrowski.
\newblock On the horizontal monotonicity of $|{\Gamma}(s)|$.
\newblock {\em Canad. Math. Bull}, 54(3):538--543, 2011.

\bibitem{2016arXiv160903682B}
R.P. {Brent}.
\newblock {On asymptotic approximations to the log-Gamma and Riemann-Siegel
  theta functions}.
\newblock {\em ArXiv e-prints}, Sept. (2016).
\newblock https://arxiv.org/abs/1609.03682v2.

\bibitem{Sai&Zven}
{Saidak} F. and P.~{Zvengrowski}.
\newblock {On the modulus of the {R}iemann {z}eta function in the critical
  strip}.
\newblock {\em Math. Slovaca}, 53:145–172, 2003.

\bibitem{Naz&Yaku}
S.~{Nazardonyavi} and S.~{Yakubovich}.
\newblock {On an inequality for the {R}iemann {Z}eta-function in the critical
  strip}.
\newblock {\em ArXiv e-prints}, June (2012).
\newblock https://arxiv.org/abs/1206.1801.

\bibitem{ZetaNotes}
M.S. Milgram.
\newblock Notes on the zeros of {R}iemann's {Z}eta function, 2011.
\newblock http://arxiv.org/abs/0911.1332.

\bibitem{Milgram}
M.S. Milgram.
\newblock Integral and series representations of {R}iemann's {Z}eta function,
  {D}irichlet's {E}ta function and a medley of related results.
\newblock {\em Journal of Mathematics, Article ID 181724}, 2013.
\newblock http://dx.doi.org/10.1155/2013/181724.

\bibitem{Titch2}
E.C. Titchmarsh and D.R Heath-Brown.
\newblock {\em The Theory of the {R}iemann {Z}eta-{F}unction}.
\newblock Oxford Science Publications, Oxford, {S}econd edition, 1986.

\bibitem{NIST:DLMF}
{NIST Digital Library of Mathematical Functions}.
\newblock http://dlmf.nist.gov/, Release 1.0.9 of 2014-08-29.

\end{thebibliography}

%\end{thebibliography}

\end{flushleft}

\end{document}